\documentclass[11pt,a4paper]{article}
\usepackage{theorem,enumerate}
\usepackage{amsmath,latexsym,amssymb,amsfonts}
\usepackage{eucal}
\usepackage{color}
\usepackage{comment}
\usepackage{cite}
\usepackage{graphics}
\theorembodyfont{\normalfont\slshape}
\newtheorem{thm}{Theorem}[section]

\newtheorem{lem}[thm]{Lemma}

\newtheorem{claim}{Claim}[section]

\newtheorem{Thm}{Theorem}
\renewcommand{\theThm}{\Alph{Thm}}

\numberwithin{equation}{section}

\newcommand{\proof}{\medbreak\noindent\textit{Proof.}\quad}
\newcommand{\Proof}[1]{\medbreak\noindent\textit{Proof of {#1}.}\quad}

\newcommand{\qed}{{$\quad\square$\vs{3.6}}}
\newcommand{\Qed}{{$\quad\blacksquare$\vs{3.6}}}

\newcommand{\Case}[1]{\medskip\noindent\textbf{Case #1:}\quad}

\newcommand{\vs}[1]{\vspace*{#1 mm}}
\newcommand{\hs}[1]{\hspace*{#1 mm}}

\def\AA{{ \mathcal{A}}}
\def\BB{{ \mathcal{B}}}
\def\CC{{ \mathcal{C}}}
\def\DD{{ \mathcal{D}}}
\def\EE{{ \mathcal{E}}}
\def\FF{{ \mathcal{F}}}
\def\GG{{ \mathcal{G}}}
\def\HH{{ \mathcal{H}}}
\def\II{{ \mathcal{I}}}
\def\JJ{{ \mathcal{J}}}
\def\KK{{ \mathcal{K}}}
\def\LL{{ \mathcal{L}}}
\def\MM{{ \mathcal{M}}}
\def\NN{{ \mathcal{N}}}
\def\OO{{ \mathcal{O}}}
\def\PP{{ \mathcal{P}}}
\def\QQ{{ \mathcal{Q}}}
\def\RR{{ \mathcal{R}}}
\def\SS{{ \mathcal{S}}}
\def\TT{{ \mathcal{T}}}
\def\UU{{ \mathcal{U}}}
\def\VV{{ \mathcal{V}}}
\def\WW{{ \mathcal{W}}}
\def\XX{{ \mathcal{X}}}
\def\YY{{ \mathcal{Y}}}
\def\ZZ{{ \mathcal{Z}}}

\renewcommand{\theenumi}{\arabic{enumi}}
\renewcommand{\labelenumi}{\upshape{(\theenumi)}}
\renewcommand{\baselinestretch}{1.3}

\bfseries\normalfont

\title{Ramsey-type results for path covers and path partitions. II. Digraphs}
\author{
Shuya Chiba$^{1}$\footnote{\texttt{e-mail:schiba@kumamoto-u.ac.jp}}\and \
Michitaka Furuya$^{2}$\footnote{\texttt{e-mail:michitaka.furuya@gmail.com}}\vs{5}\\
$^{1}$\textsl{Applied Mathematics, Faculty of Advanced Science and Technology,}\\
\textsl{Kumamoto University,}\\
\textsl{2-39-1 Kurokami, Kumamoto 860-8555, Japan}\\
$^{2}$\textsl{College of Liberal Arts and Sciences,}\\
\textsl{Kitasato University,}\\
\textsl{1-15-1 Kitasato, Minami-ku, Sagamihara, Kanagawa 252-0373, Japan}
}

\date{}

\begin{document}

\maketitle

\begin{abstract}
Recently, the authors~\cite{CF} gave Ramsey-type results for the path cover/partition number of graphs.
In this paper, we continue the research about them focusing on digraphs, and find a relationship between the path cover/partition number and forbidden structures in digraphs.

Let $D$ be a weakly connected digraph.
A family $\mathcal{P}$ of subdigraphs of $D$ is called a {\it path cover} (resp. a {\it path partition}) of $D$ if $\bigcup _{P\in \mathcal{P}}V(P)=V(D)$ (resp. $\dot\bigcup _{P\in \mathcal{P}}V(P)=V(D)$) and every element of $\mathcal{P}$ is a directed path.
The minimum cardinality of a path cover (resp. a path partition) of $D$ is denoted by  ${\rm pc}(D)$ (resp. ${\rm pp}(D)$).
In this paper, we find forbidden structure conditions assuring us that ${\rm pc}(D)$ (or ${\rm pp}(D)$) is bounded by a constant.
\end{abstract}



%%%%%%%%%%%%%%%%%%%%%%%%%%%%%%%%%%%%%%%%%%%%%%%%%%%%%%%%%%%%%%%%%%%%%%%%%%%%%%%%%%%%%%%%%%%%%%%%%%%%%%%%%%%%%%%%%%%%%%%%
%%%%%%%%%%%%%%%%%%%%%%%%%%%%%%%%%%%%%%%%%%%%%%%%%%%%%%%%%%%%%%%%%%%%%%%%%%%%%%%%%%%%%%%%%%%%%%%%%%%%%%%%%%%%%%%%%%%%%%%%
%%%%%%%%%%%%%%%%%%%%%%%%%%%%%%%%%%%%%%%%%%%%%%%%%%%%%%%%%%%%%%%%%%%%%%%%%%%%%%%%%%%%%%%%%%%%%%%%%%%%%%%%%%%%%%%%%%%%%%%%
\section{Introduction}\label{sec1}
%%%%%%%%%%%%%%%%%%%%%%%%%%%%%%%%%%%%%%%%%%%%%%%%%%%%%%%%%%%%%%%%%%%%%%%%%%%%%%%%%%%%%%%%%%%%%%%%%%%%%%%%%%%%%%%%%%%%%%%%
%%%%%%%%%%%%%%%%%%%%%%%%%%%%%%%%%%%%%%%%%%%%%%%%%%%%%%%%%%%%%%%%%%%%%%%%%%%%%%%%%%%%%%%%%%%%%%%%%%%%%%%%%%%%%%%%%%%%%%%%
%%%%%%%%%%%%%%%%%%%%%%%%%%%%%%%%%%%%%%%%%%%%%%%%%%%%%%%%%%%%%%%%%%%%%%%%%%%%%%%%%%%%%%%%%%%%%%%%%%%%%%%%%%%%%%%%%%%%%%%%

%%%%%%%%%%%%%%%%%%%%%%%%%%%%%%%%%%%%%%%%%%%%%%%%%%%%%%%%%%%%%%%%%%%%%%%%%%%%%%%%%%%%%%%%%%%%%%%%%%%%%%%%%%%%%%%%%%%%%%%%
%%%%%%%%%%%%%%%%%%%%%%%%%%%%%%%%%%%%%%%%%%%%%%%%%%%%%%%%%%%%%%%%%%%%%%%%%%%%%%%%%%%%%%%%%%%%%%%%%%%%%%%%%%%%%%%%%%%%%%%%
\subsection{Definitions}\label{sec1.1}
%%%%%%%%%%%%%%%%%%%%%%%%%%%%%%%%%%%%%%%%%%%%%%%%%%%%%%%%%%%%%%%%%%%%%%%%%%%%%%%%%%%%%%%%%%%%%%%%%%%%%%%%%%%%%%%%%%%%%%%%
%%%%%%%%%%%%%%%%%%%%%%%%%%%%%%%%%%%%%%%%%%%%%%%%%%%%%%%%%%%%%%%%%%%%%%%%%%%%%%%%%%%%%%%%%%%%%%%%%%%%%%%%%%%%%%%%%%%%%%%%

All graphs and digraphs considered in this paper are finite and simple.
Furthermore, to simplify the situations, we assume that all digraphs in this paper have no $2$-cycles (and such digraphs are sometimes called {\it oriented graphs}).
However, one might notice that our results are easily extended to digraphs with $2$-cycles.
For terms and symbols not defined in this paper, we refer the reader to \cite{D}.

Let $G$ be a graph.
Let $V(G)$ and $E(G)$ denote the {\it vertex set} and the {\it edge set} of $G$, respectively.
For a vertex $x\in V(G)$, let $N_{G}(x)$ denote the {\it neighborhood} of $x$ in $G$; thus $N_{G}(x)=\{y\in V(G): xy\in E(G)\}$.
For a subset $X$ of $V(G)$, let $N_{G}(X)=(\bigcup _{x\in X}N_{G}(x))\setminus X$.
Let $\alpha (G)$ denote the {\it independence number} of $G$, i.e., the maximum cardinality of an independent set of $G$.
Let $K_{n}$, $P_{n}$ and $K_{1,n}$ denote the {\it complete graph} of order $n$, the {\it path} of order $n$ and the {\it star} of order $n+1$, respectively.
For two positive integers $n_{1}$ and $n_{2}$, the {\it Ramsey number} $R(n_{1},n_{2})$ is the minimum positive integer $R$ such that any graph of order at least $R$ contains a clique of cardinality $n_{1}$ or an independent set of cardinality $n_{2}$.

Let $D$ be a digraph.
Let $V(D)$ and $E(D)$ denote the {\it vertex set} and the {\it arc set} of $D$, respectively.
Let $G_{D}$ be the underlying graph of $D$, i.e., $G_{D}$ is the graph on $V(D)$ with $E(G_{D})=\{xy:(x,y)\in E(D)\mbox{ or }(y,x)\in E(D)\}$.
A digraph is said to be {\it weakly connected} if its underlying graph is connected.
For a vertex $x\in V(D)$, let $d^{+}_{D}(x)$ and $d^{-}_{D}(x)$ denote the {\it out-degree} and {\it in-degree} of $x$ in $D$, respectively; thus $d^{+}_{D}(x)=|\{y\in V(D):(x,y)\in E(D)\}|$ and $d^{-}_{D}(x)=|\{y\in V(D):(y,x)\in E(D)\}|$.
A digraph $P$ is a {\it pseudo-path} if $G_{P}$ is a path.
For a pseudo-path $P$, let $r(P)=|\{x\in V(P):d_{P}^{+}(x)=2\mbox{ or }d_{P}^{-}(x)=2\}|$.
A pseudo-path $P$ is called a {\it directed path} if $r(P)=0$.
A (non-oriented) path $P$ with $V(P)=\{x_{i}:1\leq i\leq l\}$ and $E(P)=\{x_{i}x_{i+1}:1\leq i\leq l-1\}$ is often denoted by $P=x_{1}x_{2}\cdots x_{l}$.
Furthermore, if there is no fear of confusion, a directed path $P$ with $V(P)=\{x_{i}:1\leq i\leq l\}$ and $E(P)=\{(x_{i},x_{i+1}):1\leq i\leq l-1\}$ is also denoted by $P=x_{1}x_{2}\cdots x_{l}$.

For a (di)graph $\Gamma $ and a subset $X$ of $V(\Gamma )$, let $\Gamma [X]$ (resp. $\Gamma -X$) denote the sub(di)graph of $\Gamma $ induced by $X$ (resp. $V(\Gamma )\setminus X$).
For two (di)graphs $\Gamma $ and $\Gamma '$, $\Gamma $ is said to be {\it $\Gamma '$-free} if $\Gamma $ contains no induced copy of $\Gamma '$.
For a family $\GG$ of (di)graphs, a (di)graph $\Gamma $ is said to be {\it $\GG$-free} if $\Gamma $ is $\Gamma '$-free for every $\Gamma '\in \GG$.
In this context, the members of $\GG$ are called {\it forbidden sub(di)graphs}.
%For two families $\mathcal{H}_{1}$ and $\mathcal{H}_{2}$ of (di)graphs, we write $\mathcal{H}_{1}\leq \mathcal{H}_{2}$ if for every $H_{2}\in \mathcal{H}_{2}$, there exists $H_{1}\in \mathcal{H}_{1}$ such that $H_{1}$ is an induced sub(di)graph of $H_{2}$.
%The relation ``$\leq $'' between two families of forbidden subgraphs was introduced in \cite{FKLOPS}.
%Note that if $\mathcal{H}_{1}\leq \mathcal{H}_{2}$, then every $\mathcal{H}_{1}$-free graph is also $\mathcal{H}_{2}$-free.

Let $\Gamma $ be a (di)graph.
When $\Gamma $ is a graph, a family $\PP$ of subgraphs of $\Gamma $ is called a {\it path cover} of $\Gamma $ if $\bigcup _{P\in \PP}V(P)=V(\Gamma )$ and each element of $\PP$ is a path.
When $\Gamma $ is a digraph, a family $\PP$ of subdigraphs of $\Gamma $ is called a {\it path cover} of $\Gamma $ if $\bigcup _{P\in \PP}V(P)=V(\Gamma )$ and each element of $\PP$ is a directed path.
Note that some elements of a path cover of $\Gamma $ might have common vertices.
A path cover $\PP$ of $\Gamma $ is called a {\it path partition} of $\Gamma $ if the elements of $\PP$ are pairwise vertex-disjoint.
Since $\{\Gamma [\{x\}]:x\in V(\Gamma )\}$ is a path partition of $\Gamma $ (and so it is a path cover of $\Gamma $), the minimum cardinality of a path cover (or a path partition) of any (di)graph is well-defined.
The value $\min \{|\PP |:\PP\mbox{ is a path cover of }\Gamma \}$ (resp. $\min \{|\PP |:\PP\mbox{ is a path partition of }\Gamma \}$), denoted by ${\rm pc}(\Gamma )$ (resp. ${\rm pp}(\Gamma )$), is called the {\it path cover number} (resp. the {\it path partition number}) of $\Gamma $.
It is trivial that ${\rm pc}(\Gamma )\leq {\rm pp}(\Gamma )$.
Since a digraph $D$ has a Hamiltonian directed path if and only if ${\rm pp}(D)=1$, the decision problem for the path partition number is a natural generalization of the Hamiltonian directed path problem in digraphs.

%%%%%%%%%%%%%%%%%%%%%%%%%%%%%%%%%%%%%%%%%%%%%%%%%%%%%%%%%%%%%%%%%%%%%%%%%%%%%%%%%%%%%%%%%%%%%%%%%%%%%%%%%%%%%%%%%%%%%%%%
%%%%%%%%%%%%%%%%%%%%%%%%%%%%%%%%%%%%%%%%%%%%%%%%%%%%%%%%%%%%%%%%%%%%%%%%%%%%%%%%%%%%%%%%%%%%%%%%%%%%%%%%%%%%%%%%%%%%%%%%
\subsection{Ramsey-type problem and main results}\label{sec1.2}
%%%%%%%%%%%%%%%%%%%%%%%%%%%%%%%%%%%%%%%%%%%%%%%%%%%%%%%%%%%%%%%%%%%%%%%%%%%%%%%%%%%%%%%%%%%%%%%%%%%%%%%%%%%%%%%%%%%%%%%%
%%%%%%%%%%%%%%%%%%%%%%%%%%%%%%%%%%%%%%%%%%%%%%%%%%%%%%%%%%%%%%%%%%%%%%%%%%%%%%%%%%%%%%%%%%%%%%%%%%%%%%%%%%%%%%%%%%%%%%%%

The statement of a classical Ramsey theorem equals to the following:
For an integer $n\geq 1$, there exists a constant $c=c(n)$ depending on $n$ only such that every $\{K_{n},\overline{K_{n}}\}$-free graph $G$ satisfies $|V(G)|\leq c$, where $\overline{K_{n}}$ is the complement of $K_{n}$.
Its analogy for connected graphs is well-known as a folklore (see \cite[Proposition~9.4.1]{D}):
For an integer $n\geq 1$, there exists a constant $c'=c'(n)$ depending on $n$ only such that every connected $\{K_{n},K_{1,n},P_{n}\}$-free graph $G$ satisfies $|V(G)|\leq c'$.
Furthermore, some researchers have tried to extend the result to a theorem for some invariants (other than order), and found Ramsey-type properties which claim the following (see, for example, \cite{BPRS,CFKP,DDL,F}):
For a given graph-invariant $\mu $, if we impose forbidden structure conditions on a connected graph $G$, then $\mu (G)$ is bounded by a constant.

Recently, the authors~\cite{CF} gave Ramsey-type theorems for the path cover/partition number of graphs (see Theorem~\ref{ThmA} below).
In this paper, we continue the research about them focusing on digraphs, and find a relationship between the path cover/partition number and forbidden structures in digraphs.
Although
\begin{enumerate}[{$\bullet $}]
\item
the problems finding a longest directed path or a path partition in digraphs are major topics in graph theory (see \cite{BGHS,B,CMM,G,H,J}) and
\item
there are many results on forbidden subgraph conditions for the existence of a Hamiltonian path in graphs (see \cite{DGJ,FG,GH1,GH2,GH3}),
\end{enumerate}
there are only a few papers which treat directed paths in digraphs together with forbidden structure conditions.
One of its reasons is the number of forbidden structures for the existence of directed paths having a good property.
Indeed, our main theorem (Theorem~\ref{mainthm}) essentially uses eight or ten forbidden structures.
Note that half of them derives from the assumptions in Theorem~\ref{ThmA} which treats a small path cover/partition problem in graphs.
In general, when we analyze a relationship between directed paths and forbidden structures, it is necessary to previously clarify forbidden subgraph conditions for a path of graphs in quite some detail.
Thus, after getting the results in \cite{CF}, we estimate that it is a suitable time to start to analyze forbidden structure conditions for the problem finding a small path cover/partition of digraphs.

\begin{figure}
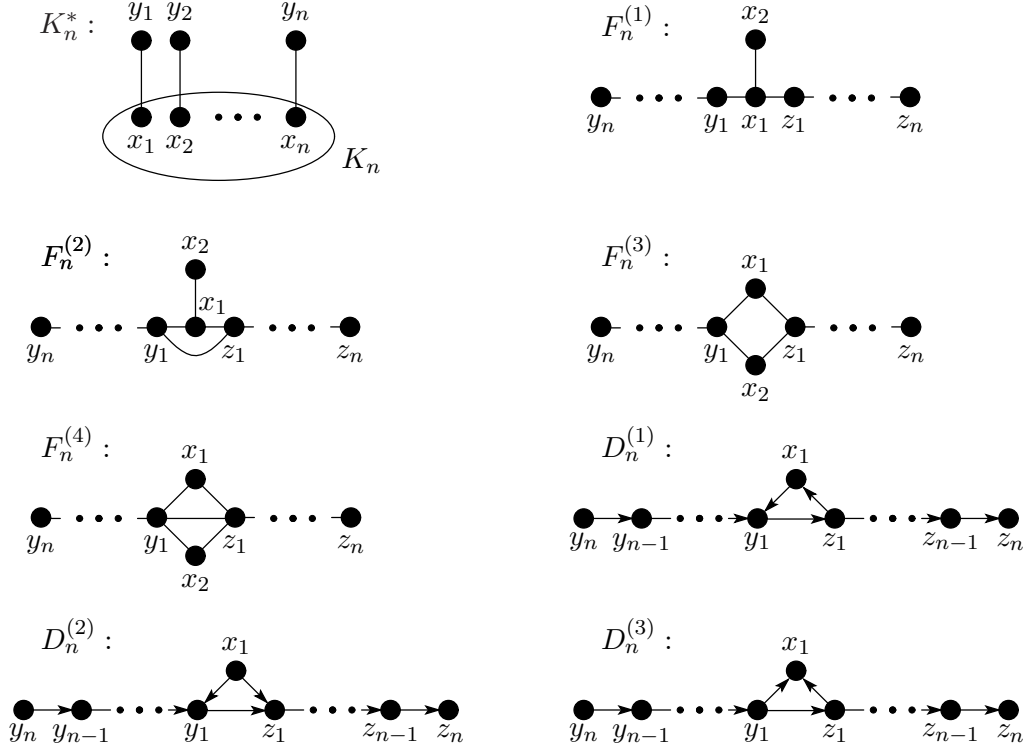

\begin{center}
%WinTpicVersion4.32a
{\unitlength 0.1in%
}%

\caption{Graphs $K^{*}_{n}$ and $F^{(i)}_{n}~(1\leq i\leq 4)$, and digraphs $D^{(i)}_{n}~(1\leq i\leq 3)$}
\label{f1}
\end{center}
\end{figure}

To state previous research and our main results, we prepare some (di)graphs which will be used as forbidden sub(di)graphs (see Figure~\ref{f1}).
Let $n$ be a positive integers.
\begin{enumerate}[{$\bullet $}]
\item
Let $K^{*}_{n}$ denote the graph with $V(K^{*}_{n})=\{x_{i},y_{i}:1\leq i\leq n\}$ and $E(K^{*}_{n})=\{x_{i}x_{j}:1\leq i<j\leq n\}\cup \{x_{i}y_{i}:1\leq i\leq n\}$.
\item
Let $A:=\{x_{1},x_{2}\}\cup \{y_{i},z_{i}:1\leq i\leq n\}$.
\begin{enumerate}[{$\circ$}]
\item
Let $F^{(1)}_{n}$ denote the graph on $A$ such that $E(F^{(1)}_{n})=\{x_{1}x_{2},x_{1}y_{1},x_{1}z_{1}\}\cup \{y_{i}y_{i+1},z_{i}z_{i+1}:1\leq i\leq n-1\}$.
\item
Let $F^{(2)}_{n}$ be the graph obtained from $F^{(1)}_{n}$ by adding the edge $y_{1}z_{1}$.
\item
Let $F^{(3)}_{n}$ denote the graph on $A$ such that $E(F^{(3)}_{n})=\{x_{1}y_{1},x_{1}z_{1},x_{2}y_{1},x_{2}z_{1}\}\cup \{y_{i}y_{i+1},z_{i}z_{i+1}:1\leq i\leq n-1\}$.
\item
Let $F^{(4)}_{n}$ be the graph obtained from $F^{(3)}_{n}$ by adding the edge $y_{1}z_{1}$.
\item
Let $D^{(1)}_{n}$ denote the digraph on $A\setminus \{x_{2}\}$ such that $E(D^{(1)}_{n})=\{(x_{1},y_{1}),(z_{1},x_{1}),(y_{1},z_{1})\}\cup \{(y_{i+1},y_{i}),(z_{i},z_{i+1}):1\leq i\leq n-1\}$.
\item
Let $D^{(2)}_{n}$ be the digraph obtained from $D^{(1)}_{n}$ by replacing the arc $(z_{1},x_{1})$ by $(x_{1},z_{1})$.
\item
Let $D^{(3)}_{n}$ be the digraph obtained from $D^{(1)}_{n}$ by replacing the arc $(x_{1},y_{1})$ by $(y_{1},x_{1})$.
\end{enumerate}
\end{enumerate}

For a fixed integer $n\geq 2$ and for a graph $G$, we consider two forbidden subgraph conditions:
\begin{enumerate}
\item[{\bf (F1)}]
$G$ is $\{K^{*}_{n},K_{1,n},F^{(1)}_{n},F^{(2)}_{n}\}$-free; and
\item[{\bf (F'1)}]
$G$ is $\{K^{*}_{n},K_{1,n},F^{(1)}_{n},F^{(2)}_{n},F^{(3)}_{n},F^{(4)}_{n}\}$-free.
\end{enumerate}
In \cite{CF}, the authors proved the following theorem (and, indeed, they also showed the forbidden subgraph conditions are necessary in a sense.).

\begin{Thm}[Chiba and Furuya~\cite{CF}]%%%%%%%%%%%%%%%%%%%%%%%%%%%%%%%%%%%%%%%%%%%%%%%%%%%%%%%%%%%%%%%%%%%%%%%%%%%%%%%%%
\label{ThmA}
For an integer $n\geq 2$, the following hold:
\begin{enumerate}
\item[{\upshape(i)}]
There exists a constant $c_{0}=c_{0}(n)$ depending on $n$ only such that ${\rm pc}(G)\leq c_{0}$ for every connected graph $G$ satisfying (F1).
\item[{\upshape(ii)}]
There exists a constant $c'_{0}=c'_{0}(n)$ depending on $n$ only such that ${\rm pp}(G)\leq c'_{0}$ for every connected graph $G$ satisfying (F'1).
\end{enumerate}
\end{Thm}
%%%%%%%%%%%%%%%%%%%%%%%%%%%%%%%%%%%%%%%%%%%%%%%%%%%%%%%%%%%%%%%%%%%%%%%%%%%%%%%%%%%%%%%%%%%%%%%%%%%%%%%%%%%%%%%%%%%%%%%%

In this paper our aim is to give an analogy of Theorem~\ref{ThmA} for digraphs.
To state our main result, we define the following four conditions for a fixed integer $n\geq 2$ and for a digraph $D$.
\begin{enumerate}
\item[{\bf (D1)}]
$G_{D}$ is $\{K^{*}_{n},K_{1,n},F^{(1)}_{n},F^{(2)}_{n}\}$-free, i.e., $G_{D}$ satisfies (F1);
\item[{\bf (D'1)}]
$G_{D}$ is $\{K^{*}_{n},K_{1,n},F^{(1)}_{n},F^{(2)}_{n},F^{(3)}_{n},F^{(4)}_{n}\}$-free, i.e., $G_{D}$ satisfies (F'1);
\item[{\bf (D2)}]
$D$ is $\{D^{(1)}_{n},D^{(2)}_{n},D^{(3)}_{n}\}$-free; and
\item[{\bf (D3)}]
$r(P)\leq n$ for every induced pseudo-path $P$ of $D$.
\end{enumerate}
Our main theorem is the following.

\begin{thm}%%%%%%%%%%%%%%%%%%%%%%%%%%%%%%%%%%%%%%%%%%%%%%%%%%%%%%%%%%%%%%%%%%%%%%%%%%%%%%%%%%%%%%%%%%%%%%%%%%%%%%%%%%%%%
\label{mainthm}
For an integer $n\geq 2$, the following hold:
\begin{enumerate}
\item[{\upshape(i)}]
There exists a constant $c_{1}=c_{1}(n)$ depending on $n$ only such that ${\rm pc}(D)\leq c_{1}$ for every weakly connected digraph $D$ satisfying (D1), (D2) and (D3).
\item[{\upshape(ii)}]
There exists a constant $c'_{1}=c'_{1}(n)$ depending on $n$ only such that ${\rm pp}(D)\leq c'_{1}$ for every weakly connected digraph $D$ satisfying (D'1), (D2) and (D3).
\end{enumerate}
\end{thm}
%%%%%%%%%%%%%%%%%%%%%%%%%%%%%%%%%%%%%%%%%%%%%%%%%%%%%%%%%%%%%%%%%%%%%%%%%%%%%%%%%%%%%%%%%%%%%%%%%%%%%%%%%%%%%%%%%%%%%%%%

Strictly speaking, the condition (D3) is not a forbidden subdigraph condition.
However, when we consider conditions for bounding the path cover/partition number of digraphs, we must forbid a pseudo-path $P$ with large $r(P)$ because ${\rm pc}(P)={\rm pp}(P)=r(P)+1$ for every pseudo-path $P$.
Since no finite forbidden subdigraph conditions can exclude all pseudo-paths $P$ with large $r(P)$, the condition (D3) is reasonable in our main results.

The condition that ``weakly connected digraph $D$ satisfying (D1) (or (D'1)), (D2) and (D3)'' is best possible in a sense because essentially weaker conditions than the conditions (D1) (resp. (D'1)), (D2) and (D3) cannot assure us the existence of a small path cover (resp. a small path partition) of digraphs.
(We omit its details, but the discussion in \cite{CF} will be a good reference for the checking.)

We introduce some lemmas in Section~\ref{sec2}, and prove Theorem~\ref{mainthm} in Section~\ref{sec3}.
We remark that the assumptions imposed on digraphs in Theorem~\ref{mainthm} are closely divided into two parts:
\begin{enumerate}[$\bullet $]
\item
Fundamental assumptions for the underlying graph $G_{D}$ ((D1) and (D'1)), which lead to the existence of small path cover/partitions of $G_{D}$, and
\item
additional assumptions which are peculiar to digraphs ((D2) and (D3)).
\end{enumerate}
However in general, when we compare Ramsey-type problems on graphs and digraphs, the above situation does not always occur.
To demonstrate it, in Section~\ref{sec4}, we give a similar theorem on the cycle cover/partition number of digraphs and show that the assumptions for the underlying graph do not assure us the existence of small cycle cover/partitions of graphs.

%%%%%%%%%%%%%%%%%%%%%%%%%%%%%%%%%%%%%%%%%%%%%%%%%%%%%%%%%%%%%%%%%%%%%%%%%%%%%%%%%%%%%%%%%%%%%%%%%%%%%%%%%%%%%%%%%%%%%%%%
%%%%%%%%%%%%%%%%%%%%%%%%%%%%%%%%%%%%%%%%%%%%%%%%%%%%%%%%%%%%%%%%%%%%%%%%%%%%%%%%%%%%%%%%%%%%%%%%%%%%%%%%%%%%%%%%%%%%%%%%
%%%%%%%%%%%%%%%%%%%%%%%%%%%%%%%%%%%%%%%%%%%%%%%%%%%%%%%%%%%%%%%%%%%%%%%%%%%%%%%%%%%%%%%%%%%%%%%%%%%%%%%%%%%%%%%%%%%%%%%%
\section{Lemmas}\label{sec2}
%%%%%%%%%%%%%%%%%%%%%%%%%%%%%%%%%%%%%%%%%%%%%%%%%%%%%%%%%%%%%%%%%%%%%%%%%%%%%%%%%%%%%%%%%%%%%%%%%%%%%%%%%%%%%%%%%%%%%%%%
%%%%%%%%%%%%%%%%%%%%%%%%%%%%%%%%%%%%%%%%%%%%%%%%%%%%%%%%%%%%%%%%%%%%%%%%%%%%%%%%%%%%%%%%%%%%%%%%%%%%%%%%%%%%%%%%%%%%%%%%
%%%%%%%%%%%%%%%%%%%%%%%%%%%%%%%%%%%%%%%%%%%%%%%%%%%%%%%%%%%%%%%%%%%%%%%%%%%%%%%%%%%%%%%%%%%%%%%%%%%%%%%%%%%%%%%%%%%%%%%%

%%%%%%%%%%%%%%%%%%%%%%%%%%%%%%%%%%%%%%%%%%%%%%%%%%%%%%%%%%%%%%%%%%%%%%%%%%%%%%%%%%%%%%%%%%%%%%%%%%%%%%%%%%%%%%%%%%%%%%%%
%%%%%%%%%%%%%%%%%%%%%%%%%%%%%%%%%%%%%%%%%%%%%%%%%%%%%%%%%%%%%%%%%%%%%%%%%%%%%%%%%%%%%%%%%%%%%%%%%%%%%%%%%%%%%%%%%%%%%%%%
\subsection{Global structures}\label{sec2-1}
%%%%%%%%%%%%%%%%%%%%%%%%%%%%%%%%%%%%%%%%%%%%%%%%%%%%%%%%%%%%%%%%%%%%%%%%%%%%%%%%%%%%%%%%%%%%%%%%%%%%%%%%%%%%%%%%%%%%%%%%
%%%%%%%%%%%%%%%%%%%%%%%%%%%%%%%%%%%%%%%%%%%%%%%%%%%%%%%%%%%%%%%%%%%%%%%%%%%%%%%%%%%%%%%%%%%%%%%%%%%%%%%%%%%%%%%%%%%%%%%%

In this subsection, we introduce some known lemmas.

As we mentioned in Subsection~\ref{sec1.1}, we have ${\rm pc}(D)\leq {\rm pp}(D)$ for a digraph $D$.
This, together with a result by Gallai and Milgram~\cite{GM}, leads to the following lemma.

\begin{lem}%%%%%%%%%%%%%%%%%%%%%%%%%%%%%%%%%%%%%%%%%%%%%%%%%%%%%%%%%%%%%%%%%%%%%%%%%%%%%%%%%%%%%%%%%%%%%%%%%%%%%%%%%%%%%
\label{lem2.0}
For a digraph $D$, ${\rm pc}(D)\leq {\rm pp}(D)\leq \alpha (G_{D})$.
\end{lem}
%%%%%%%%%%%%%%%%%%%%%%%%%%%%%%%%%%%%%%%%%%%%%%%%%%%%%%%%%%%%%%%%%%%%%%%%%%%%%%%%%%%%%%%%%%%%%%%%%%%%%%%%%%%%%%%%%%%%%%%%

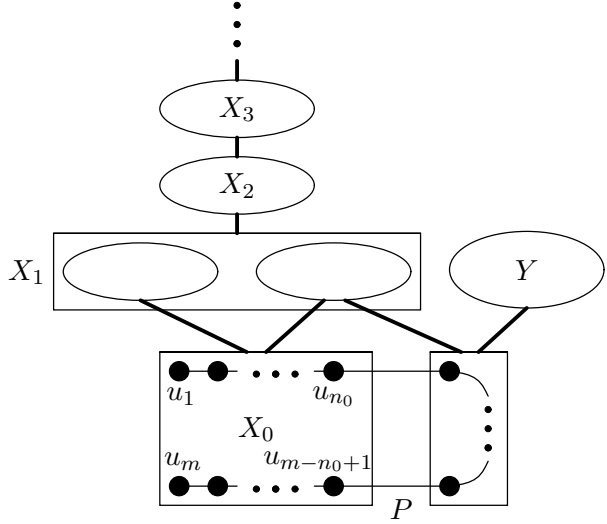
\begin{figure}
\begin{center}
%WinTpicVersion4.32a
{\unitlength 0.1in%
\begin{picture}(31.9000,26.1000)(2.1000,-29.0000)%
% CIRCLE 2 0 0 0 Black Black  
% 4 1200 2200 1200 2250 1200 2250 1200 2250
% 
\special{sh 1.000}%
\special{ia 1200 2200 50 50 0.0000000 6.2831853}%
\special{pn 8}%
\special{ar 1200 2200 50 50 0.0000000 6.2831853}%
% CIRCLE 2 0 0 0 Black Black  
% 4 1400 2200 1400 2250 1400 2250 1400 2250
% 
\special{sh 1.000}%
\special{ia 1400 2200 50 50 0.0000000 6.2831853}%
\special{pn 8}%
\special{ar 1400 2200 50 50 0.0000000 6.2831853}%
% CIRCLE 2 0 0 0 Black Black  
% 4 2000 2200 2000 2250 2000 2250 2000 2250
% 
\special{sh 1.000}%
\special{ia 2000 2200 50 50 0.0000000 6.2831853}%
\special{pn 8}%
\special{ar 2000 2200 50 50 0.0000000 6.2831853}%
% CIRCLE 2 0 0 0 Black Black  
% 4 2600 2200 2600 2250 2600 2250 2600 2250
% 
\special{sh 1.000}%
\special{ia 2600 2200 50 50 0.0000000 6.2831853}%
\special{pn 8}%
\special{ar 2600 2200 50 50 0.0000000 6.2831853}%
% CIRCLE 2 0 0 0 Black Black  
% 4 2600 2800 2600 2850 2600 2850 2600 2850
% 
\special{sh 1.000}%
\special{ia 2600 2800 50 50 0.0000000 6.2831853}%
\special{pn 8}%
\special{ar 2600 2800 50 50 0.0000000 6.2831853}%
% CIRCLE 2 0 0 0 Black Black  
% 4 2000 2800 2000 2850 2000 2850 2000 2850
% 
\special{sh 1.000}%
\special{ia 2000 2800 50 50 0.0000000 6.2831853}%
\special{pn 8}%
\special{ar 2000 2800 50 50 0.0000000 6.2831853}%
% CIRCLE 2 0 0 0 Black Black  
% 4 1400 2800 1400 2850 1400 2850 1400 2850
% 
\special{sh 1.000}%
\special{ia 1400 2800 50 50 0.0000000 6.2831853}%
\special{pn 8}%
\special{ar 1400 2800 50 50 0.0000000 6.2831853}%
% CIRCLE 2 0 0 0 Black Black  
% 4 1200 2800 1200 2850 1200 2850 1200 2850
% 
\special{sh 1.000}%
\special{ia 1200 2800 50 50 0.0000000 6.2831853}%
\special{pn 8}%
\special{ar 1200 2800 50 50 0.0000000 6.2831853}%
% BOX 2 0 3 0 Black Black  
% 2 1100 2100 2200 2900
% 
\special{pn 8}%
\special{pa 1100 2100}%
\special{pa 2200 2100}%
\special{pa 2200 2900}%
\special{pa 1100 2900}%
\special{pa 1100 2100}%
\special{pa 2200 2100}%
\special{fp}%
% LINE 2 0 3 0 Black Black  
% 4 1200 2200 1500 2200 1900 2200 2000 2200
% 
\special{pn 8}%
\special{pa 1200 2200}%
\special{pa 1500 2200}%
\special{fp}%
\special{pa 1900 2200}%
\special{pa 2000 2200}%
\special{fp}%
% LINE 2 0 3 0 Black Black  
% 4 1200 2800 1500 2800 1900 2800 2000 2800
% 
\special{pn 8}%
\special{pa 1200 2800}%
\special{pa 1500 2800}%
\special{fp}%
\special{pa 1900 2800}%
\special{pa 2000 2800}%
\special{fp}%
% DOT 0 0 3 0 Black Black  
% 4 1600 2815 1800 2815 1700 2815 1700 2815
% 
\special{pn 4}%
\special{sh 1}%
\special{ar 1600 2815 16 16 0 6.2831853}%
\special{sh 1}%
\special{ar 1800 2815 16 16 0 6.2831853}%
\special{sh 1}%
\special{ar 1700 2815 16 16 0 6.2831853}%
\special{sh 1}%
\special{ar 1700 2815 16 16 0 6.2831853}%
% DOT 0 0 3 0 Black Black  
% 4 1600 2215 1800 2215 1700 2215 1700 2215
% 
\special{pn 4}%
\special{sh 1}%
\special{ar 1600 2215 16 16 0 6.2831853}%
\special{sh 1}%
\special{ar 1800 2215 16 16 0 6.2831853}%
\special{sh 1}%
\special{ar 1700 2215 16 16 0 6.2831853}%
\special{sh 1}%
\special{ar 1700 2215 16 16 0 6.2831853}%
% DOT 0 0 3 0 Black Black  
% 4 2800 2500 2800 2400 2800 2600 2800 2600
% 
\special{pn 4}%
\special{sh 1}%
\special{ar 2800 2500 16 16 0 6.2831853}%
\special{sh 1}%
\special{ar 2800 2400 16 16 0 6.2831853}%
\special{sh 1}%
\special{ar 2800 2600 16 16 0 6.2831853}%
\special{sh 1}%
\special{ar 2800 2600 16 16 0 6.2831853}%
% SPLINE 2 0 3 0 Black Black  
% 4 2600 2200 2750 2250 2800 2330 2800 2330
% 
\special{pn 8}%
\special{pa 2600 2200}%
\special{pa 2634 2205}%
\special{pa 2668 2211}%
\special{pa 2699 2220}%
\special{pa 2727 2233}%
\special{pa 2751 2251}%
\special{pa 2771 2275}%
\special{pa 2787 2303}%
\special{pa 2800 2330}%
\special{fp}%
% SPLINE 2 0 3 0 Black Black  
% 4 2600 2800 2750 2750 2800 2670 2800 2670
% 
\special{pn 8}%
\special{pa 2600 2800}%
\special{pa 2634 2795}%
\special{pa 2668 2789}%
\special{pa 2699 2780}%
\special{pa 2727 2767}%
\special{pa 2751 2749}%
\special{pa 2771 2725}%
\special{pa 2787 2697}%
\special{pa 2800 2670}%
\special{fp}%
% STR 2 0 3 0 Black Black  
% 4 1600 2400 1600 2500 5 0 0 0
% $X_{0}$
\put(16.0000,-25.0000){\makebox(0,0){$X_{0}$}}%
% BOX 2 0 3 0 Black Black  
% 2 2500 2100 2900 2900
% 
\special{pn 8}%
\special{pa 2500 2100}%
\special{pa 2900 2100}%
\special{pa 2900 2900}%
\special{pa 2500 2900}%
\special{pa 2500 2100}%
\special{pa 2900 2100}%
\special{fp}%
% LINE 2 0 3 0 Black Black  
% 4 2000 2200 2600 2200 2600 2800 2000 2800
% 
\special{pn 8}%
\special{pa 2000 2200}%
\special{pa 2600 2200}%
\special{fp}%
\special{pa 2600 2800}%
\special{pa 2000 2800}%
\special{fp}%
% STR 2 0 3 0 Black Black  
% 4 2350 2820 2350 2920 5 0 0 0
% $P$
\put(23.5000,-29.2000){\makebox(0,0){$P$}}%
% STR 2 0 3 0 Black Black  
% 4 1210 2230 1210 2330 5 0 0 0
% $u_{1}$
\put(12.1000,-23.3000){\makebox(0,0){$u_{1}$}}%
% STR 2 0 3 0 Black Black  
% 4 2000 2230 2000 2330 5 0 0 0
% $u_{n_{0}}$
\put(20.0000,-23.3000){\makebox(0,0){$u_{n_{0}}$}}%
% STR 2 0 3 0 Black Black  
% 4 1920 2570 1920 2670 5 0 0 0
% $u_{m-n_{0}+1}$
\put(19.2000,-26.7000){\makebox(0,0){$u_{m-n_{0}+1}$}}%
% STR 2 0 3 0 Black Black  
% 4 1220 2570 1220 2670 5 0 0 0
% $u_{m}$
\put(12.2000,-26.7000){\makebox(0,0){$u_{m}$}}%
% ELLIPSE 2 0 3 0 Black Black  
% 4 3000 1670 3400 1870 3400 1870 3400 1870
% 
\special{pn 8}%
\special{ar 3000 1670 400 200 0.0000000 6.2831853}%
% LINE 0 0 3 0 Black Black  
% 2 2750 2100 3000 1870
% 
\special{pn 20}%
\special{pa 2750 2100}%
\special{pa 3000 1870}%
\special{fp}%
% ELLIPSE 2 0 3 0 Black Black  
% 4 1500 1230 1900 1380 1900 1380 1900 1380
% 
\special{pn 8}%
\special{ar 1500 1230 400 150 0.0000000 6.2831853}%
% ELLIPSE 2 0 3 0 Black Black  
% 4 1500 830 1900 980 1900 980 1900 980
% 
\special{pn 8}%
\special{ar 1500 830 400 150 0.0000000 6.2831853}%
% STR 2 0 3 0 Black Black  
% 4 1500 1130 1500 1230 5 0 0 0
% $X_{2}$
\put(15.0000,-12.3000){\makebox(0,0){$X_{2}$}}%
% STR 2 0 3 0 Black Black  
% 4 1500 730 1500 830 5 0 0 0
% $X_{3}$
\put(15.0000,-8.3000){\makebox(0,0){$X_{3}$}}%
% LINE 0 0 3 0 Black Black  
% 2 1500 1080 1500 980
% 
\special{pn 20}%
\special{pa 1500 1080}%
\special{pa 1500 980}%
\special{fp}%
% LINE 0 0 3 0 Black Black  
% 2 1500 680 1500 580
% 
\special{pn 20}%
\special{pa 1500 680}%
\special{pa 1500 580}%
\special{fp}%
% DOT 0 0 3 0 Black Black  
% 4 1500 390 1500 290 1500 490 1500 490
% 
\special{pn 4}%
\special{sh 1}%
\special{ar 1500 390 16 16 0 6.2831853}%
\special{sh 1}%
\special{ar 1500 290 16 16 0 6.2831853}%
\special{sh 1}%
\special{ar 1500 490 16 16 0 6.2831853}%
\special{sh 1}%
\special{ar 1500 490 16 16 0 6.2831853}%
% LINE 0 0 3 0 Black Black  
% 2 1500 1480 1500 1380
% 
\special{pn 20}%
\special{pa 1500 1480}%
\special{pa 1500 1380}%
\special{fp}%
% STR 2 0 3 0 Black Black  
% 4 410 1580 410 1680 5 0 0 0
% $X_{1}$
\put(4.1000,-16.8000){\makebox(0,0){$X_{1}$}}%
% STR 2 0 3 0 Black Black  
% 4 3000 1570 3000 1670 5 0 0 0
% $Y$
\put(30.0000,-16.7000){\makebox(0,0){$Y$}}%
% ELLIPSE 2 0 3 0 Black Black  
% 4 1000 1680 1400 1830 1400 1830 1400 1830
% 
\special{pn 8}%
\special{ar 1000 1680 400 150 0.0000000 6.2831853}%
% ELLIPSE 2 0 3 0 Black Black  
% 4 2000 1680 2400 1830 2400 1830 2400 1830
% 
\special{pn 8}%
\special{ar 2000 1680 400 150 0.0000000 6.2831853}%
% LINE 0 0 3 0 Black Black  
% 6 2660 2100 2060 1830 1950 1830 1650 2100 1550 2100 1000 1830
% 
\special{pn 20}%
\special{pa 2660 2100}%
\special{pa 2060 1830}%
\special{fp}%
\special{pa 1950 1830}%
\special{pa 1650 2100}%
\special{fp}%
\special{pa 1550 2100}%
\special{pa 1000 1830}%
\special{fp}%
% BOX 2 0 3 0 Black Black  
% 2 550 1480 2450 1880
% 
\special{pn 8}%
\special{pa 550 1480}%
\special{pa 2450 1480}%
\special{pa 2450 1880}%
\special{pa 550 1880}%
\special{pa 550 1480}%
\special{pa 2450 1480}%
\special{fp}%
\end{picture}}%

\caption{Path $P$, and sets $X_{i}$ and $Y$}
\label{f2}
\end{center}
\end{figure}

Now we pick up two useful lemmas proved in \cite{CF}.
In the remainder of this subsection, we fix an integer $n\geq 3$ and a connected $\{K^{*}_{n},K_{1,n},F^{(1)}_{n},F^{(2)}_{n}\}$-free graph $G$.
Set $n_{0}=\max\{\lceil \frac{n^{2}-n-2}{2}\rceil ,n\}$.
Take a longest induced path $P$ of $G$, and write $P=u_{1}u_{2}\cdots u_{m}$.
Let $X_{0}=\{u_{i}:1\leq i\leq n_{0}\mbox{ or }m-n_{0}+1\leq i\leq m\}$ and $Y=N_{G}(V(P)\setminus X_{0})\setminus (X_{0}\cup N_{G}(X_{0}))$.
Note that if $|V(P)|\leq 2n_{0}$, then $X_{0}=V(P)$ and $Y=\emptyset $.
We further remark that $N_{G}(y)\cap V(P)\subseteq \{u_{i}:n_{0}+1\leq i\leq m-n_{0}\}$ for every $y\in Y$.
%For each $i$ with $n_{0}+1\leq i\leq m-n_{0}$, let $Y_{i}=\{y\in Y:\min\{j:n_{0}+1\leq j\leq m-n_{0},~yu_{j}\in E(G)\}=i\}$.
Now we recursively define the sets $X_{i}~(i\geq 1)$ as follows:
Let $X_{1}=N_{G}(X_{0})\setminus V(P)$, and for $i$ with $i\geq 2$, let $X_{i}=N_{G}(X_{i-1})\setminus (V(P)\cup Y\cup (\bigcup _{1\leq j\leq i-2}X_{j}))$ (see Figure~\ref{f2}).
Then $N_{G}(V(P))$ is the disjoint union of $X_{1}$ and $Y$, and the following lemma holds.

\begin{lem}[Chiba and Furuya~{\cite[Lemma~2.6]{CF}}]%%%%%%%%%%%%%%%%%%%%%%%%%%%%%%%%%%%%%%%%%%%%%%%%%%%%%%%%%%%%%%%%%%%%
\label{lem2.4}
The set $V(G)$ is the disjoint union of $V(P)$, $Y$ and $X_{i}~(1\leq i\leq 2n_{0}-1)$.
In particular, $X_{i}=\emptyset $ for all $i$ with $i\geq 2n_{0}$.
\end{lem}
%%%%%%%%%%%%%%%%%%%%%%%%%%%%%%%%%%%%%%%%%%%%%%%%%%%%%%%%%%%%%%%%%%%%%%%%%%%%%%%%%%%%%%%%%%%%%%%%%%%%%%%%%%%%%%%%%%%%%%%%

We recursively define the values $\alpha _{i}~(i\geq 0)$ as follows:
Let $\alpha _{0}=2\lceil \frac{n_{0}}{2} \rceil $, and for $i$ with $i\geq 1$, let $\alpha _{i}=(n-1)R(n,\alpha _{i-1}+1)-1$.
We remark that $\alpha _{i}$ depends on $n$ and $i$.
The following lemma holds.

\begin{lem}[Chiba and Furuya~{\cite[Lemma~2.7]{CF}}]%%%%%%%%%%%%%%%%%%%%%%%%%%%%%%%%%%%%%%%%%%%%%%%%%%%%%%%%%%%%%%%%%%%%
\label{lem-alpha}
For an integer $i$ with $i\geq 0$, $\alpha (G[X_{i}])\leq \alpha _{i}$.
\end{lem}
%%%%%%%%%%%%%%%%%%%%%%%%%%%%%%%%%%%%%%%%%%%%%%%%%%%%%%%%%%%%%%%%%%%%%%%%%%%%%%%%%%%%%%%%%%%%%%%%%%%%%%%%%%%%%%%%%%%%%%%%

%%%%%%%%%%%%%%%%%%%%%%%%%%%%%%%%%%%%%%%%%%%%%%%%%%%%%%%%%%%%%%%%%%%%%%%%%%%%%%%%%%%%%%%%%%%%%%%%%%%%%%%%%%%%%%%%%%%%%%%%
%%%%%%%%%%%%%%%%%%%%%%%%%%%%%%%%%%%%%%%%%%%%%%%%%%%%%%%%%%%%%%%%%%%%%%%%%%%%%%%%%%%%%%%%%%%%%%%%%%%%%%%%%%%%%%%%%%%%%%%%
\subsection{An induced path and its neighbors}\label{sec2-2}
%%%%%%%%%%%%%%%%%%%%%%%%%%%%%%%%%%%%%%%%%%%%%%%%%%%%%%%%%%%%%%%%%%%%%%%%%%%%%%%%%%%%%%%%%%%%%%%%%%%%%%%%%%%%%%%%%%%%%%%%
%%%%%%%%%%%%%%%%%%%%%%%%%%%%%%%%%%%%%%%%%%%%%%%%%%%%%%%%%%%%%%%%%%%%%%%%%%%%%%%%%%%%%%%%%%%%%%%%%%%%%%%%%%%%%%%%%%%%%%%%

In this subsection, we fix an integer $n\geq 3$ and an $\{F^{(1)}_{n},F^{(2)}_{n}\}$-free graph $G$.
Take an induced path $Q$ of $G$, and write $Q=v_{1}v_{2}\cdots v_{l}$.
For a vertex $y\in N_{G}(V(Q))$, let $i_{y}=\min\{i:1\leq i\leq l,~yv_{i}\in E(G)\}$ and $j_{y}=\max\{i:1\leq i\leq l,~yv_{i}\in E(G)\}$.
Since $N_{G}(y)\cap V(Q)\neq \emptyset $, the indices $i_{y}$ and $j_{y}$ are well-defined.

\begin{lem}%%%%%%%%%%%%%%%%%%%%%%%%%%%%%%%%%%%%%%%%%%%%%%%%%%%%%%%%%%%%%%%%%%%%%%%%%%%%%%%%%%%%%%%%%%%%%%%%%%%%%%%%%%%%%
\label{lem-2-3-0}
Let $y\in N_{G}(V(Q))$ be a vertex such that $N_{G}(y)\cap V(Q)\subseteq \{v_{i}:n+1\leq i\leq l-n\}$.
Then the following hold:
\begin{enumerate}
\item[{\upshape(i)}]
We have $1\leq j_{y}-i_{y}\leq 3$.
\item[{\upshape(ii)}]
If $G$ is $F^{(3)}_{n}$-free, then $N_{G}(y)\cap V(Q)=\{v_{i}:i_{y}\leq i\leq j_{y}\}$.
\end{enumerate}
\end{lem}
%%%%%%%%%%%%%%%%%%%%%%%%%%%%%%%%%%%%%%%%%%%%%%%%%%%%%%%%%%%%%%%%%%%%%%%%%%%%%%%%%%%%%%%%%%%%%%%%%%%%%%%%%%%%%%%%%%%%%%%%
\proof
If $j_{y}-i_{y}\geq 3$ and $\{yv_{i_{y}+1},yv_{j_{y}-1}\}\not\subseteq E(G)$, say $yv_{i_{y}+1}\notin E(G)$, then
$$
\{v_{i_{y}-n},v_{i_{y}-n+1},\ldots ,v_{i_{y}+1},y,v_{j_{y}},v_{j_{y}+1},\ldots ,v_{j_{y}+n-2}\}
$$
induces a copy of $F^{(1)}_{n}$ in $G$, which is a contradiction.
Thus
\begin{align}
\mbox{if $j_{y}-i_{y}\geq 3$, then $yv_{i_{y}+1},yv_{j_{y}-1}\in E(G)$.}\label{cond-lem-2-3-0-1}
\end{align}
\begin{enumerate}
\item[{\upshape(i)}]
If $i_{y}=j_{y}$, then $\{v_{i_{y}-n},v_{i_{y}-n+1},\ldots ,v_{i_{y}},y,v_{i_{y}+1},\ldots ,v_{i_{y}+n}\}$ induces a copy of $F^{(1)}_{n}$ in $G$, which is a contradiction.
Thus $j_{y}-i_{y}\geq 1$.
Suppose that $j_{y}-i_{y}\geq 4$.
Then by (\ref{cond-lem-2-3-0-1}), we have $v_{i_{y}+1}y\in E(G)$.
If $yv_{i_{y}+2}\in E(G)$, then $\{v_{i_{y}-n+1},v_{i_{y}-n+2},\ldots ,v_{i_{y}},y,v_{i_{y}+2},v_{j_{y}},v_{j_{y}+1},\ldots ,v_{j_{y}+n-1}\}$ induces a copy of $F^{(1)}_{n}$ in $G$; if $yv_{i_{y}+2}\notin E(G)$, then $\{v_{i_{y}-n+1},v_{i_{y}-n+2},\ldots ,v_{i_{y}+2},y,v_{j_{y}},v_{j_{y}+1},\ldots ,v_{j_{y}+n-2}\}$ induces a copy of $F^{(2)}_{n}$ in $G$.
In either case, we obtain a contradiction.
This implies that $j_{y}-i_{y}\leq 3$.
\item[{\upshape(ii)}]
Suppose that $G$ is $F^{(3)}_{n}$-free and $N_{G}(y)\cap V(Q)\neq \{v_{i}:i_{y}\leq i\leq j_{y}\}$.
By the definitions of $i_{y}$ and $j_{y}$, $N_{G}(y)\cap V(Q)\subseteq \{v_{i}:i_{y}\leq i\leq j_{y}\}$ and $yv_{i_{y}},yv_{j_{y}}\in E(G)$.
Thus there exists an index $i$ with $i_{y}+1\leq i\leq j_{y}-1$ such that $yv_{i}\notin E(G)$.
Considering (\ref{cond-lem-2-3-0-1}) and (i), this forces $j_{y}-i_{y}=2$ and $yv_{i_{y}+1}\notin E(G)$.
Then $\{v_{i_{y}-n+1},v_{i_{y}-n+2},\ldots ,v_{i_{y}},y,v_{i_{y}+1},v_{i_{y}+2},\ldots ,v_{i_{y}+n+1}\}$ induces a copy of $F^{(3)}_{n}$ in $G$, which is a contradiction.
\qed
\end{enumerate}

\begin{lem}%%%%%%%%%%%%%%%%%%%%%%%%%%%%%%%%%%%%%%%%%%%%%%%%%%%%%%%%%%%%%%%%%%%%%%%%%%%%%%%%%%%%%%%%%%%%%%%%%%%%%%%%%%%%%
\label{lem-2-3-00}
Let $y,y'\in N_{G}(V(Q))$ be vertices with $y\neq y'$ and $(N_{G}(y)\cup N_{G}(y'))\cap V(Q)\subseteq \{v_{i}:n+1\leq i\leq l-n\}$.
Then the following hold:
\begin{enumerate}
\item[{\upshape(i)}]
If $i_{y'}\geq i_{y}+n+3$, then $yy'\notin E(G)$.
\item[{\upshape(ii)}]
If $G$ is $\{F^{(3)}_{n},F^{(4)}_{n}\}$-free and $i_{y}=i_{y'}$, then $yy'\in E(G)$.
\end{enumerate}
\end{lem}
%%%%%%%%%%%%%%%%%%%%%%%%%%%%%%%%%%%%%%%%%%%%%%%%%%%%%%%%%%%%%%%%%%%%%%%%%%%%%%%%%%%%%%%%%%%%%%%%%%%%%%%%%%%%%%%%%%%%%%%%
\proof
\begin{enumerate}
\item[{\upshape(i)}]
Suppose that $i_{y'}\geq i_{y}+n+3$ and $yy'\in E(G)$.
By Lemma~\ref{lem-2-3-0}(i), $j_{y}\leq i_{y}+3\leq (i_{y'}-n-3)+3=i_{y'}-n$, and hence
\begin{align}
N_{G}(y)\cap \{v_{i}:i_{y'}-n+1\leq i\leq l\}=\emptyset .\label{cond-lem-2-3-00-1}
\end{align}
Again by Lemma~\ref{lem-2-3-0}(i), we have $1\leq j_{y'}-i_{y'}\leq 3$.
Considering (\ref{cond-lem-2-3-00-1}), if $j_{y'}-i_{y'}=1$, then
$$
\{v_{i_{y'}-n+1},v_{i_{y'}-n+2},\ldots ,v_{i_{y'}},y',y,v_{i_{y'}+1},v_{i_{y'}+2},\ldots ,v_{i_{y'}+n}\}
$$
induces a copy of $F^{(2)}_{n}$ in $G$; if $2\leq j_{y'}-i_{y'}\leq 3$, then
$$
\{v_{i_{y'}-n+1},v_{i_{y'}-n+2},\ldots ,v_{i_{y'}},y',y,v_{j_{y'}},v_{j_{y'}+1},\ldots ,v_{j_{y'}+n-1}\}
$$
induces a copy of $F^{(1)}_{n}$ in $G$.
In either case, we obtain a contradiction.

\item[{\upshape(ii)}]
Suppose that $G$ is $\{F^{(3)}_{n},F^{(4)}_{n}\}$-free, $i_{y}=i_{y'}$ and $yy'\notin E(G)$.
Without loss of generality, we may assume that $j_{y}\geq j_{y'}$.
It follows from Lemma~\ref{lem-2-3-0}, we have $yv_{i_{y}+1},y'v_{i_{y}+1}\in E(G)$.
If $j_{y}=j_{y'}=i_{y}+1$, then
$$
\{v_{i_{y}-n+1},v_{i_{y}-n+2},\ldots ,v_{i_{y}},y,y',v_{i_{y}+1},v_{i_{y}+2},\ldots ,v_{i_{y}+n}\}
$$
induces a copy of $F^{(4)}_{n}$ in $G$; if $j_{y}=j_{y'}\geq i_{y}+2$, then
$$
\{v_{i_{y}-n+1},v_{i_{y}-n+2},\ldots ,v_{i_{y}},y,y',v_{j_{y}},v_{j_{y}+1},\ldots ,v_{j_{y}+n-1}\}
$$
induces a copy of $F^{(3)}_{n}$ in $G$; if $j_{y}>j_{y'}$, then
$$
\{v_{i_{y}-n},v_{i_{y}-n+1},\ldots ,v_{i_{y}},y',y,v_{j_{y}},v_{j_{y}+1},\ldots ,v_{j_{y}+n-2}\}
$$
induces a copy of $F^{(1)}_{n}$ in $G$.
In any case, we obtain a contradiction.
\qed
\end{enumerate}

%%%%%%%%%%%%%%%%%%%%%%%%%%%%%%%%%%%%%%%%%%%%%%%%%%%%%%%%%%%%%%%%%%%%%%%%%%%%%%%%%%%%%%%%%%%%%%%%%%%%%%%%%%%%%%%%%%%%%%%%
%%%%%%%%%%%%%%%%%%%%%%%%%%%%%%%%%%%%%%%%%%%%%%%%%%%%%%%%%%%%%%%%%%%%%%%%%%%%%%%%%%%%%%%%%%%%%%%%%%%%%%%%%%%%%%%%%%%%%%%%
\subsection{An induced directed path and its neighbors}\label{sec2-3}
%%%%%%%%%%%%%%%%%%%%%%%%%%%%%%%%%%%%%%%%%%%%%%%%%%%%%%%%%%%%%%%%%%%%%%%%%%%%%%%%%%%%%%%%%%%%%%%%%%%%%%%%%%%%%%%%%%%%%%%%
%%%%%%%%%%%%%%%%%%%%%%%%%%%%%%%%%%%%%%%%%%%%%%%%%%%%%%%%%%%%%%%%%%%%%%%%%%%%%%%%%%%%%%%%%%%%%%%%%%%%%%%%%%%%%%%%%%%%%%%%

In this subsection, we prove the following lemma which plays a key role in the proof of Theorem~\ref{mainthm}.

\begin{lem}%%%%%%%%%%%%%%%%%%%%%%%%%%%%%%%%%%%%%%%%%%%%%%%%%%%%%%%%%%%%%%%%%%%%%%%%%%%%%%%%%%%%%%%%%%%%%%%%%%%%%%%%%%%%%
\label{lem-2-3-1}
Let $n\geq 3$ be an integer, and let $D$ be a digraph satisfying (D1), (D2) and (D3).
Let $Q$ be an induced directed path in $D$, and write $Q=v_{1}v_{2}\cdots v_{l}$.
Let $Y'$ be a subset of $N_{G_{D}}(V(Q))$ such that $N_{G_{D}}(y)\cap V(Q)\subseteq \{v_{i}:n+1\leq i\leq l-n\}$ for every $y\in Y'$.
Then the following hold:
\begin{enumerate}
\item[{\upshape(i)}]
We have ${\rm pc}(D[V(Q)\cup Y'])\leq \frac{(n-2)n(n+5)}{2}+6(n-2)$.
\item[{\upshape(ii)}]
If $D$ satisfies (D'1), then ${\rm pp}(D[V(Q)\cup Y'])\leq \frac{n(n+5)}{2}+1$ and ${\rm pp}(D[(V(Q)\setminus \{v_{t}\})\cup Y'])\leq \frac{n(n+5)}{2}+1$ for each $t\in \{1,l\}$.
\end{enumerate}
\end{lem}
%%%%%%%%%%%%%%%%%%%%%%%%%%%%%%%%%%%%%%%%%%%%%%%%%%%%%%%%%%%%%%%%%%%%%%%%%%%%%%%%%%%%%%%%%%%%%%%%%%%%%%%%%%%%%%%%%%%%%%%%
\proof
Since $Q=v_{1}v_{2}\cdots v_{l}$ is a directed path, $(v_{i},v_{i+1})\in E(D)$ for every $i$ with $1\leq i\leq l-1$.
If $Y'=\emptyset $, then
\begin{enumerate}[{$\bullet $}]
\item
$\{Q\}$ is a path partition of $D[V(Q)\cup Y']$ and
\item
$\{Q-\{v_{t}\}\}$ is a path partition of $D[(V(Q)\setminus \{v_{t}\})\cup Y']$ for each $t\in \{1,l\}$,
\end{enumerate}
as desired.
Thus we may assume that $Y'\neq \emptyset $.

For each $y\in Y'$, let $i_{y}=\min\{i:n+1\leq i\leq l-n,~yv_{i}\in E(G_{D})\}$ and $j_{y}=\max\{i:n+1\leq i\leq l-n,~yv_{i}\in E(G_{D})\}$.
Since $N_{G_{D}}(y)\cap \{v_{i}:n+1\leq i\leq l-n\}\neq \emptyset $, the indices $i_{y}$ and $j_{y}$ are well-defined.
Since $G_{D}$ is $\{F^{(1)}_{n},F^{(2)}_{n}\}$-free, it follows from Lemma~\ref{lem-2-3-0}(i) that $1\leq j_{y}-i_{y}\leq 3$.
Furthermore, if $D$ satisfies (D'1), then by Lemma~\ref{lem-2-3-0}(ii), we have $N_{G_{D}}(y)\cap V(Q)=\{v_{i}:i_{y}\leq i\leq j_{y}\}$.

\begin{figure}
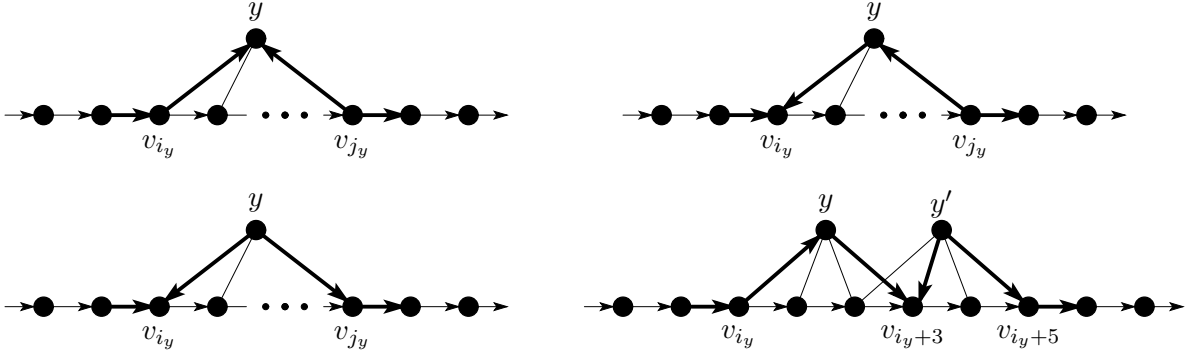

\begin{center}
%WinTpicVersion4.32a
{\unitlength 0.1in%
}%

\caption{Bad sets and a related induced pseudo-path in the proof of Claim~\ref{cl-lem-2-3-2}}
\label{f3}
\end{center}
\end{figure}

We call a subset $B$ of $Y'$ {\it bad} if either
\begin{enumerate}[{$\bullet $}]
\item
$B=\{y\}$ for some $y\in Y'$, and $E(D)\cap \{(y,v_{i_{y}}),(v_{j_{y}},y)\}\neq \emptyset $; or
\item
$B=\{y,y'\}$ for some $y,y'\in Y'$ with $y\neq y'$, and the following conditions hold:
$yy'\notin E(G_{D})$, $j_{y}-i_{y}=j_{y'}-i_{y'}=3$, $i_{y'}=i_{y}+2$, and $\{(v_{i_{y}},y),(y,v_{i_{y}+3}),(y',v_{i_{y}+3}),(y',v_{i_{y}+5})\}\subseteq E(D)$
\end{enumerate}
(see Figure~\ref{f3}).
If $\{y\}$ is a bad set and $j_{y}-i_{y}=1$, then
$$
\{v_{i_{y}-n+1},v_{i_{y}-n+2},\ldots ,v_{i_{y}},y,v_{i_{y}+1},v_{i_{y}+2},\ldots ,v_{i_{y}+n}\}
$$
induces a copy of one of $D^{(1)}_{n}$, $D^{(2)}_{n}$ and $D^{(3)}_{n}$ in $D$, which contradicts the assumption that $D$ satisfies (D2).
Thus if $\{y\}$ is a bad set, then $2\leq j_{y}-i_{y}\leq 3$.

Let $i_{B}=\min\{i_{y}:y\in B\}$ for each bad set $B$, and let $I=\{i_{B}:B\mbox{ is a bad set}\}$.

\begin{claim}%%%%%%%%%%%%%%%%%%%%%%%%%%%%%%%%%%%%%%%%%%%%%%%%%%%%%%%%%%%%%%%%%%%%%%%%%%%%%%%%%%%%%%%%%%%%%%%%%%%%%%%%%%%
\label{cl-lem-2-3-2}
We have $|I|\leq \frac{n(n+5)}{2}$.
\end{claim}
%%%%%%%%%%%%%%%%%%%%%%%%%%%%%%%%%%%%%%%%%%%%%%%%%%%%%%%%%%%%%%%%%%%%%%%%%%%%%%%%%%%%%%%%%%%%%%%%%%%%%%%%%%%%%%%%%%%%%%%%
\proof
We may assume that $I\neq \emptyset $.
Now we consider an operation recursively defining the indices $h_{1},h_{2},\ldots $ with $h_{1}<h_{2}<\cdots $ as follows:
Let $h_{1}=\min\{i:i\in I\}$.
For $p\geq 2$, we assume that the index $h_{p-1}$ has been defined.
If there exists an index $i\in I$ with $i\geq h_{p-1}+n+5$, we let $h_{p}=\min\{i\in I:i\geq h_{p-1}+n+5\}$; otherwise, we finish the operation.

Let $I'=\{h_{p}:p\geq 1\}$, and set $p_{0}=|I'|$.
For each $p$ with $1\leq p\leq p_{0}$, we take a bad set $B_{p}$ with $i_{B_{p}}=h_{p}$, and let $j_{p}=\max\{j_{y}:y\in B_{p}\}$.
Then for $y\in B_{p}$,
\begin{align}
\begin{cases}
i_{y}=i_{B_{p}} & (|B_{p}|=1)\\
i_{y}\in \{i_{B_{p}},i_{B_{p}}+2\} & (|B_{p}|=2),
\end{cases}
\label{cond-cl-lem-2-3-2-bad01}
\end{align}
and
\begin{align}
\begin{cases}
j_{p}\leq h_{p}+3 & (|B_{p}|=1)\\
j_{p}=h_{p}+5 & (|B_{p}|=2).
\end{cases}
\label{cond-cl-lem-2-3-2-bad02}
\end{align}

Fix two indices $p$ and $p'$ with $1\leq p<p'\leq p_{0}$, and let $y\in B_{p}$ and $y'\in B_{p'}$.
Then by (\ref{cond-cl-lem-2-3-2-bad01}), we have
$$
i_{y}+n+3\leq (i_{B_{p}}+2)+n+3=h_{p}+n+5\leq h_{p'}=i_{B_{p'}}\leq i_{y'}.
$$
This together with Lemma~\ref{lem-2-3-00}(i) implies that $yy'\notin E(G_{D})$.
Consequently,
\begin{align}
\mbox{there is no edge of $G_{D}$ between $B_{p}$ and $B_{p'}$ for $1\leq p<p'\leq p_{0}$.}\label{cond-cl-lem-2-3-2-1}
\end{align}

Suppose that $p_{0}\geq \frac{n+1}{2}$.
For each $p$ with $1\leq p\leq p_{0}$, we define a pseudo-path $R_{p}$ in $D$ as follows:
If $|B_{p}|=1$, say $B_{p}=\{y\}$, let $R_{p}=v_{i_{y}-1}v_{i_{y}}yv_{j_{y}}v_{j_{y}+1}$; if $|B_{p}|=2$, say $B_{p}=\{y,y'\}$ with $i_{y'}=i_{y}+2$, let $R_{p}=v_{i_{y}-1}v_{i_{y}}yv_{i_{y}+3}y'v_{i_{y}+5}v_{i_{y}+6}$ (see Figure~\ref{f3} again).
Then $R_{p}$ is an induced pseudo-path in $D$ and $r(R_{p})=2$.
Let $R$ be the subdigraph of $D$ induced by
$$
(V(Q)\setminus \{v_{i}:1\leq p\leq p_{0},~h_{p}-1\leq i\leq j_{p}+1\})\cup \left(\bigcup _{1\leq p\leq p_{0}}V(R_{p})\right).
$$
By (\ref{cond-cl-lem-2-3-2-1}) and the definition of a bad set, $\bigcup _{1\leq p\leq p_{0}}B_{p}$ is an independent set of $G_{D}$.
For $p$ and $p'$ with $1\leq p<p'\leq p_{0}$, it follows from (\ref{cond-cl-lem-2-3-2-bad02}) that $j_{p}+1\leq h_{p}+6\leq (h_{p'}-n-5)+6=h_{p'}-n+1\leq h_{p'}-2$, and hence $V(R_{p})\cap V(R_{p'})=\emptyset $.
Therefore, $R$ is an induced pseudo-path in $D$ and $r(R)=\sum _{1\leq p\leq p_{0}}r(R_{p})=2p_{0}\geq n+1$, which contradicts the assumption that $D$ satisfies (D3).
Thus $p_{0}\leq \frac{n}{2}$.

Suppose that $I\setminus \{i:1\leq p\leq p_{0},~h_{p}\leq i\leq h_{p}+n+4\}\neq \emptyset $, and take its minimum element $h$.
Since $h\neq h_{1}~(=\min\{i:i\in I\})$ and $h\in I$, we have $h_{1}<h$.
In particular, the value $p^{*}=\max\{p:1\leq p\leq p_{0},~h_{p}<h\}$ is well-defined.
Then by the definition of $h$,
\begin{align}
h\neq h_{p^{*}+1}\label{cond-hp*-000}
\end{align}
and $h\notin \{i:h_{p^{*}}\leq i\leq h_{p^{*}}+n+4\}$, and so
\begin{align}
h\geq h_{p^{*}}+n+5.\label{cond-hp*-1}
\end{align}
By the maximality of $p^{*}$,
\begin{align}
\mbox{no index $i$ with $h_{p^{*}}+1\leq i\leq h-1$ equals to $h_{p^{*}+1}$.}\label{cond-hp*-2}
\end{align}
Considering the definition of the indices $h_{p}~(1\leq p\leq p_{0})$, it follows from (\ref{cond-hp*-1}), (\ref{cond-hp*-2}) and the minimality of $h$ that $h_{p^{*}+1}=h$, which contradicts (\ref{cond-hp*-000}).
Thus $I\subseteq \{i:1\leq p\leq p_{0},~h_{p}\leq i\leq h_{p}+n+4\}$, and so $|I|\leq p_{0}(n+5)\leq \frac{n(n+5)}{2}$.
\qed

For each $i$ with $n+1\leq i\leq l-n-1$, let $Y'_{i}=\{y\in Y':i_{y}=i\}$.
Recall that
\begin{align}
1\leq j_{y}-i_{y}\leq 3~~\mbox{and}~~n+1\leq i_{y}<j_{y}\leq l-n ~~~\mbox{for every }y\in Y'.\label{cond-fund-lem22-index}
\end{align}
By (\ref{cond-fund-lem22-index}), $\{y\in Y':i_{y}\leq n\mbox{ and }i_{y}\geq l-n\}=\emptyset $, and hence $Y'$ is the disjoint union of $Y'_{i}~(n+1\leq i\leq l-n-1)$.
Note that $Y'_{i}\subseteq N_{G_{D}}(v_{i})$ and $Y'_{i}\cap N_{G_{D}}(v_{i-1})=\emptyset $.
Hence if $\alpha (G_{D}[Y'_{i}])\geq n-1$, then $G_{D}[\{v_{i-1},v_{i}\}\cup Y'_{i}]$ contains a copy of $K_{1,n}$ as an induced subgraph, which contradicts the $K_{1,n}$-freeness of $G_{D}$.
Thus
\begin{align}
\alpha (G_{D}[Y'_{i}])\leq n-2\mbox{ for every }i\mbox{ with }n+1\leq i\leq l-n-1.\label{cond-lem-2-3-1-0000}
\end{align}
Furthermore, it follows from Lemma~\ref{lem-2-3-00}(ii) that
\begin{align}
\mbox{if $D$ satisfies (D'1), then $Y'_{i}$ is a clique of $G_{D}$ for every $i$ with $n+1\leq i\leq l-n-1$.}\label{cond-lem-2-3-1-0001}
\end{align}
Let $\tilde{Y}=\bigcup _{i\in I}Y'_{i}$.
By the definition of $I$, $\tilde{Y}\cap B\neq \emptyset $ for every bad set $B$, and so
\begin{align}
B\not\subseteq Y'\setminus \tilde{Y} \mbox{ for every bad set }B.\label{cond-lem-2-3-1-00}
\end{align}

\begin{figure}
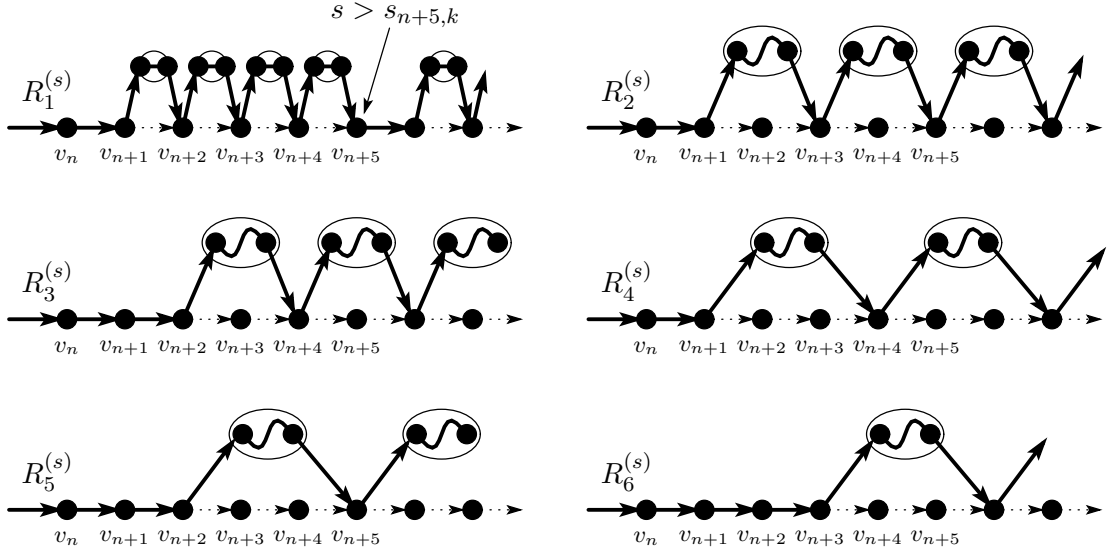

\begin{center}
%WinTpicVersion4.32a
{\unitlength 0.1in%
}%

\caption{The directed paths $R^{(s)}_{a}$ }
\label{f5}
\end{center}
\end{figure}

Now we prove that ${\rm pc}(D[V(Q)\cup Y'])\leq \frac{(n-2)n(n+5)}{2}+6(n-2)$, which proves (i).
By Lemma~\ref{lem2.0}, Claim~\ref{cl-lem-2-3-2} and (\ref{cond-lem-2-3-1-0000}), we have ${\rm pc}(D[\tilde{Y}])\leq \alpha (G_{D}[\tilde{Y}])\leq \sum _{i\in I}\alpha (G_{D}[Y'_{i}])\leq \frac{(n-2)n(n+5)}{2}$.
Thus, to complete the proof of (i), it suffices to show that ${\rm pc}(D[V(Q)\cup (Y'\setminus \tilde{Y})])\leq 6(n-2)$.
For two integers $i$ and $k$ with  $n+1\leq i\leq l-n-k$ and $k\in \{1,2,3\}$, we define a set $Y'_{i,k}$ as follows:
If $i\notin I$, let $Y'_{i,k}=\{y\in Y'_{i}:j_{y}=i_{y}+k\}$; if $i\in I$, let $Y'_{i,k}=\emptyset $.
%Since there is no edge of $G_{D}$ between $Y$ and $\{v_{i}:l-n+1\leq i\leq l\}$, $Y'_{i,j}=\emptyset $ if $i+j\geq l-n+1$.
%Since there is no edge of $G_{D}$ between $Y$ and $\{v_{i}:l-n+1\leq i\leq l\}$, it follows from Lemma~\ref{lem-2-3-0}(i) implies that $Y'\setminus \tilde{Y}$ is the disjoint union of $Y'_{i,j}~(n+1\leq i\leq l-n-j,~j\in \{1,2,3\})$.
Then by (\ref{cond-fund-lem22-index}), $Y'\setminus \tilde{Y}$ is the disjoint union of $Y'_{i,k}~(n+1\leq i\leq l-n-k,~k\in \{1,2,3\})$.
By Lemma~\ref{lem2.0} and (\ref{cond-lem-2-3-1-0000}), there exists a path cover $\QQ_{i,k}=\{Q^{(1)}_{i,k},Q^{(2)}_{i,k},\ldots ,Q^{(s_{i,k})}_{i,k}\}$ of $D[Y'_{i,k}]$ with $s_{i,k}\leq n-2$, where $\QQ_{i,k}=\emptyset $ and $s_{i,k}=0$ if $Y'_{i,k}=\emptyset $.
Let $p^{(s)}_{i,k}$ and $q^{(s)}_{i,k}$ be the internal vertex and terminal vertex of $Q^{(s)}_{i,k}$, respectively.
If $Y'_{i,k}\neq \emptyset $, then for every $y\in Y'_{i,k}$, $\{y\}$ is not a bad set, and so $(v_{i},y),(y,v_{i+k})\in E(D)$.
Hence for each integer $s$ with $1\leq s\leq n-2$, we can define a directed path $R^{(s)}_{i,k}$ as follows:
If $s\leq s_{i,k}$, let $R^{(s)}_{i,k}=v_{i}p^{(s)}_{i,k}Q^{(s)}_{i,k}q^{(s)}_{i,k}v_{i+k}$; otherwise, let $R^{(s)}_{i,k}=v_{i}v_{i+1}\cdots v_{i+k}$.
For $r\in \{0,1\}$, we define the value $\xi _{r}$ as
$$
\xi _{r}=
\begin{cases}
l-n-r & (l\mbox{ is odd})\\
l-n-1+r & (l\mbox{ is even}).
\end{cases}
$$
We define the values $\xi'_{0}$, $\xi'_{1}$ and $\xi'_{2}$ as
$$
\xi '_{0}=
\begin{cases}
l-n-2 & (l-2n\equiv 0~(\mbox{mod }3))\\
l-n & (l-2n\equiv 1~(\mbox{mod }3))\\
l-n-1 & (l-2n\equiv 2~(\mbox{mod }3)),
\end{cases}
$$
$$
\xi '_{1}=
\begin{cases}
l-n-1 & (l-2n\equiv 0~(\mbox{mod }3))\\
l-n-2 & (l-2n\equiv 1~(\mbox{mod }3))\\
l-n & (l-2n\equiv 2~(\mbox{mod }3)),
\end{cases}~~~
$$
and
$$
\xi '_{2}=
\begin{cases}
l-n & (l-2n\equiv 0~(\mbox{mod }3))\\
l-n-1 & (l-2n\equiv 1~(\mbox{mod }3))\\
l-n-2 & (l-2n\equiv 2~(\mbox{mod }3)).
\end{cases}~~~
$$
Let
\begin{align*}
R^{(s)}_{1} &= v_{1}v_{2}\cdots v_{n+1}R^{(s)}_{n+1,1}v_{n+2}R^{(s)}_{n+2,1}v_{n+3}\cdots v_{l-n-1}R^{(s)}_{l-n-1,1}v_{l-n}v_{l-n+1}\cdots v_{l},\\
R^{(s)}_{2} &= v_{1}v_{2}\cdots v_{n+1}R^{(s)}_{n+1,2}v_{n+3}R^{(s)}_{n+3,2}v_{n+5}\cdots v_{\xi _{0}-2}R^{(s)}_{\xi _{0}-2,2}v_{\xi _{0}}v_{\xi _{0}+1}\cdots v_{l},\\
R^{(s)}_{3} &= v_{1}v_{2}\cdots v_{n+2}R^{(s)}_{n+2,2}v_{n+4}R^{(s)}_{n+4,2}v_{n+6}\cdots v_{\xi _{1}-2}R^{(s)}_{\xi _{1}-2,2}v_{\xi _{1}}v_{\xi _{1}+1}\cdots v_{l},\\
R^{(s)}_{4} &= v_{1}v_{2}\cdots v_{n+1}R^{(s)}_{n+1,3}v_{n+4}R^{(s)}_{n+4,3}v_{n+7}\cdots v_{\xi '_{0}-3}R^{(s)}_{\xi '_{0}-3,3}v_{\xi '_{0}}v_{\xi '_{0}+1}\cdots v_{l},\\
R^{(s)}_{5} &= v_{1}v_{2}\cdots v_{n+2}R^{(s)}_{n+2,3}v_{n+5}R^{(s)}_{n+5,3}v_{n+8}\cdots v_{\xi '_{1}-3}R^{(s)}_{\xi '_{1}-3,3}v_{\xi '_{1}}v_{\xi '_{1}+1}\cdots v_{l},\mbox{ and}\\
R^{(s)}_{6} &= v_{1}v_{2}\cdots v_{n+3}R^{(s)}_{n+3,3}v_{n+6}R^{(s)}_{n+6,3}v_{n+9}\cdots v_{\xi '_{2}-3}R^{(s)}_{\xi '_{2}-3,3}v_{\xi '_{2}}v_{\xi '_{2}+1}\cdots v_{l}
\end{align*}
(see Figure~\ref{f5}).
Then we easily verify that $\{R^{(s)}_{a}:1\leq a\leq 6,~1\leq s\leq n-2\}$ is a path cover of $D[V(Q)\cup (Y'\setminus \tilde{Y})]$ having cardinality at most $6(n-2)$.
Consequently, ${\rm pc}(D[V(Q)\cup (Y'\setminus \tilde{Y})])\leq 6(n-2)$, which proves (i).

Next we assume that $D$ satisfies (D'1), and prove that (ii) holds.
By Lemma~\ref{lem2.0}, Claim~\ref{cl-lem-2-3-2} and (\ref{cond-lem-2-3-1-0001}), we have ${\rm pp}(D[\tilde{Y}])\leq \alpha (G_{D}[\tilde{Y}])\leq \sum _{i\in I}\alpha (G_{D}[Y'_{i}])\leq |I|\leq \frac{n(n+5)}{2}$.
Thus, to complete the proof of (ii), it suffices to show that
\begin{align}
{\rm pp}(D[V(Q)\cup (Y'\setminus \tilde{Y})])=1\label{cond-dilem-pp1}
\end{align}
and
\begin{align}
{\rm pp}(D[(V(Q)\setminus \{v_{t}\})\cup (Y'\setminus \tilde{Y})])=1 \mbox{ for each }t\in \{1,l\}.\label{cond-dilem-pp2}
\end{align}

\begin{figure}
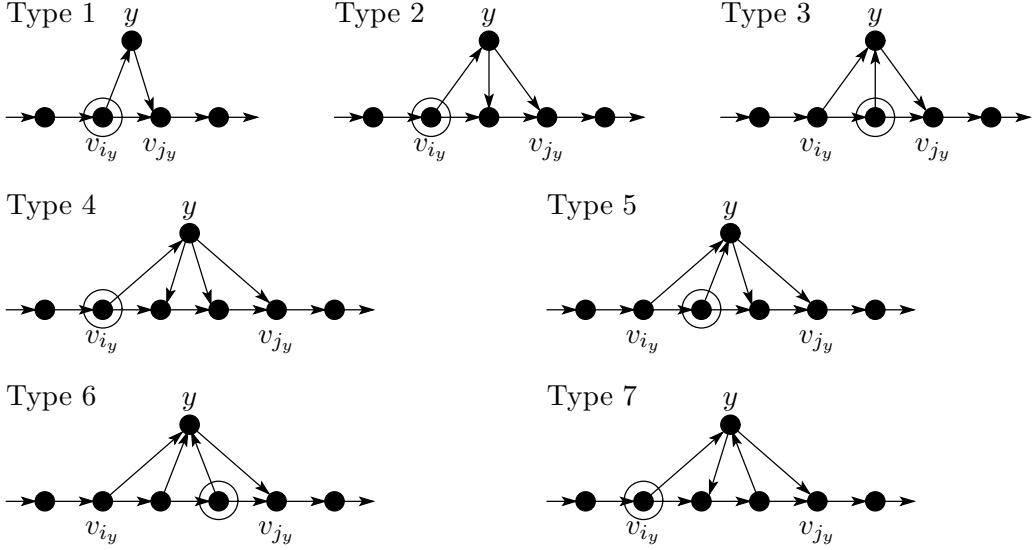

\begin{center}
%WinTpicVersion4.32a
{\unitlength 0.1in%
}%

\caption{Seven types of $y\in Y'\setminus \tilde{Y}$ where the vertex $v_{\beta _{y}}$ is circled.}
\label{f4}
\end{center}
\end{figure}

Now we fix a vertex $y\in Y'\setminus \tilde{Y}$.
Recall that $1\leq j_{y}-i_{y}\leq 3$ and $N_{G_{D}}(y)\cap V(Q)=\{v_{i}:i_{y}\leq i\leq j_{y}\}$.
Since $\{y\}$ is not a bad set, $(v_{i_{y}},y),(y,v_{j_{y}})\in E(D)$.
In particular, the index
$$
\beta _{y}=\min\{i:i_{y}\leq i\leq j_{y}-1,~(v_{i},y),(y,v_{i+1})\in E(D)\}
$$
is well-defined.
Furthermore, each vertex $y$ in $Y'\setminus \tilde{Y}$ is divided into the following seven types (see Figure~\ref{f4}):\\
{\bf Type 1}~~$j_{y}=i_{y}+1$.\\
{\bf Type 2}~~$j_{y}=i_{y}+2$ and $(y,v_{i_{y}+1})\in E(D)$.\\
{\bf Type 3}~~$j_{y}=i_{y}+2$ and $(v_{i_{y}+1},y)\in E(D)$.\\
{\bf Type 4}~~$j_{y}=i_{y}+3$ and $(y,v_{i_{y}+1}),(y,v_{i_{y}+2})\in E(D)$.\\
{\bf Type 5}~~$j_{y}=i_{y}+3$ and $(v_{i_{y}+1},y),(y,v_{i_{y}+2})\in E(D)$.\\
{\bf Type 6}~~$j_{y}=i_{y}+3$ and $(v_{i_{y}+1},y),(v_{i_{y}+2},y)\in E(D)$.\\
{\bf Type 7}~~$j_{y}=i_{y}+3$ and $(y,v_{i_{y}+1}),(v_{i_{y}+2},y)\in E(D)$.\\
Note that $\beta _{y}=i_{y}$ if the type of $y$ is 1, 2, 4 or 7; $\beta _{y}=i_{y}+1$ if the type of $y$ is 3 or 5; $\beta _{y}=i_{y}+2$ if the type of $y$ is 6.
For $i$ with $n+1\leq i\leq l-n-1$, let $W_{i}=\{y\in Y'\setminus \tilde{Y}:\beta _{y}=i\}$.
Then $Y'\setminus \tilde{Y}$ is the disjoint union of $W_{i}~(n+1\leq i\leq l-n-1)$.

\begin{claim}%%%%%%%%%%%%%%%%%%%%%%%%%%%%%%%%%%%%%%%%%%%%%%%%%%%%%%%%%%%%%%%%%%%%%%%%%%%%%%%%%%%%%%%%%%%%%%%%%%%%%%%%%%%
\label{cl-lem-2-3-3}
For $i$ with $n+1\leq i\leq l-n-1$, $W_{i}$ is a clique of $G_{D}$.
\end{claim}
%%%%%%%%%%%%%%%%%%%%%%%%%%%%%%%%%%%%%%%%%%%%%%%%%%%%%%%%%%%%%%%%%%%%%%%%%%%%%%%%%%%%%%%%%%%%%%%%%%%%%%%%%%%%%%%%%%%%%%%%
\proof
By way of contradiction, suppose that $W_{i}$ is not a clique of $G_{D}$.
Then there exist two vertices $y$ and $y'$ in $W_{i}$ with $yy'\notin E(G_{D})$ whose types are $t_{y}$ and $t_{y'}$, respectively.
We may assume that $t_{y}\leq t_{y'}$.
By Lemma~\ref{lem-2-3-00}(ii), $i_{y}\neq i_{y'}$.
This implies that $t_{y}\neq t_{y'}$ and
$$
(t_{y},t_{y'})\notin \{(1,2),(1,4),(1,7),(2,4),(2,7),(4,7),(3,5)\}.
$$
Thus there are $\binom{7}{2}-7=14$ estimated cases on the pair $(t_{y},t_{y'})$.
We divide the proof into five cases.

\medskip
\noindent
\textbf{Case 1:} $(t_{y},t_{y'})\in \{(1,3),(1,6),(2,5),(3,6)\}$.

Note that $i_{y}>i_{y'}$ and $j_{y}=j_{y'}$.
Hence
$$
\{v_{i_{y'}-n+2},v_{i_{y'}-n+3},\ldots ,v_{i_{y'}},y',v_{j_{y}},y,v_{j_{y}+1},\ldots ,v_{j_{y}+n}\}
$$
induces a copy of $F^{(1)}_{n}$ in $G_{D}$, which is a contradiction.

\medskip
\noindent
\textbf{Case 2:} $(t_{y},t_{y'})=(1,5)$.

The set
$$
\{v_{i_{y'}-n+1},v_{i_{y'}-n+2},\ldots ,v_{i_{y'}+1},y,y',v_{j_{y'}},v_{j_{y'}+1},\ldots ,v_{j_{y'}+n-2}\}
$$
induces a copy of $F^{(2)}_{n}$ in $G_{D}$, which is a contradiction.

\medskip
\noindent
\textbf{Case 3:} $(t_{y},t_{y'})\in \{(2,3),(2,6),(4,5),(5,6)\}$.

Note that $i_{y}>i_{y'}$, $j_{y'}-i_{y'}\geq 2$, $j_{y'}=j_{y}-1$ and $(y,v_{j_{y}-1})\in E(D)$.
Hence
$$
\{v_{i_{y'}-n+3},v_{i_{y'}-n+4},\ldots ,v_{i_{y'}},y',v_{j_{y}-1},y,v_{j_{y}},v_{j_{y}+1},\ldots ,v_{j_{y}+n-1}\}
$$
induces a copy of $D^{(2)}_{n}$ in $D$, which is a contradiction.

\medskip
\noindent
\textbf{Case 4:} $(t_{y},t_{y'})\in \{(3,4),(3,7),(5,7)\}$.

Note that $i_{y}=i_{y'}-1$, $(v_{i_{y}+1},y)\in E(D)$ and $j_{y'}-i_{y'}=3$ (and so $j_{y}<j_{y'}$).
Hence
$$
\{v_{i_{y}-n+1},v_{i_{y}-n+2},\ldots ,v_{i_{y}},y,v_{i_{y}+1},y',v_{j_{y'}},v_{j_{y'}+1},\ldots ,v_{j_{y'}+n-3}\}
$$
induces a copy of $D^{(3)}_{n}$ in $D$, which is a contradiction.

\medskip
\noindent
\textbf{Case 5:} $(t_{y},t_{y'})\in \{(4,6),(6,7)\}$.

It follows that $\{y,y'\}$ is a bad set and $\{y,y'\}\subseteq W_{i}\subseteq Y'\setminus \tilde{Y}$, which contradicts (\ref{cond-lem-2-3-1-00}).

This completes the proof of Claim~\ref{cl-lem-2-3-3}.
\qed

Fix an index $i$ with $n+1\leq i\leq l-n-1$, and we define the directed path $R_{i}$ from $v_{i}$ to $v_{i+1}$ with $V(R_{i})=\{v_{i},v_{i+1}\}\cup W_{i}$ as follows:
If $W_{i}=\emptyset $, let $R_{i}$ be the directed path $v_{i}v_{i+1}$.
We assume that $W_{i}\neq \emptyset $.
By Lemma~\ref{lem2.0} and Claim~\ref{cl-lem-2-3-3}, $D$ contains a directed path $R'_{i}$ with $V(R'_{i})=W_{i}$.
Let $R_{i}=v_{i}p_{i}R'_{i}q_{i}v_{i+1}$ where $p_{i}$ and $q_{i}$ are is the internal vertex and terminal vertex of $R'_{i}$, respectively.
Since $(v_{i},y),(y,v_{i+1})\in E(D)$ for every $y\in W_{i}$, $R_{i}$ is a directed path in $D$.
Hence
$$
R^{*}=v_{1}v_{2}\cdots v_{n+1}R_{n+1}v_{n+2}R_{n+2}v_{n+3}\cdots v_{l-n-1}R_{l-n-1}v_{l-n}v_{l-n+1}\cdots v_{l}
$$
is a directed path in $D$ with $V(R^{*})=V(Q)\cup (\bigcup _{n+1\leq i\leq l-n-1}V(R_{i}))~(=V(Q)\cup (Y'\setminus \tilde{Y}))$.
Furthermore, $R^{*}-v_{1}$ and $R^{*}-v_{l}$ are also directed paths in $D$.
Consequently, (\ref{cond-dilem-pp1}) and (\ref{cond-dilem-pp2}) hold.

This completes the proof of Lemma~\ref{lem-2-3-1}.
\qed

%%%%%%%%%%%%%%%%%%%%%%%%%%%%%%%%%%%%%%%%%%%%%%%%%%%%%%%%%%%%%%%%%%%%%%%%%%%%%%%%%%%%%%%%%%%%%%%%%%%%%%%%%%%%%%%%%%%%%%%%
%%%%%%%%%%%%%%%%%%%%%%%%%%%%%%%%%%%%%%%%%%%%%%%%%%%%%%%%%%%%%%%%%%%%%%%%%%%%%%%%%%%%%%%%%%%%%%%%%%%%%%%%%%%%%%%%%%%%%%%%
%%%%%%%%%%%%%%%%%%%%%%%%%%%%%%%%%%%%%%%%%%%%%%%%%%%%%%%%%%%%%%%%%%%%%%%%%%%%%%%%%%%%%%%%%%%%%%%%%%%%%%%%%%%%%%%%%%%%%%%%
\section{Proof of Theorem~\ref{mainthm}}\label{sec3}
%%%%%%%%%%%%%%%%%%%%%%%%%%%%%%%%%%%%%%%%%%%%%%%%%%%%%%%%%%%%%%%%%%%%%%%%%%%%%%%%%%%%%%%%%%%%%%%%%%%%%%%%%%%%%%%%%%%%%%%%
%%%%%%%%%%%%%%%%%%%%%%%%%%%%%%%%%%%%%%%%%%%%%%%%%%%%%%%%%%%%%%%%%%%%%%%%%%%%%%%%%%%%%%%%%%%%%%%%%%%%%%%%%%%%%%%%%%%%%%%%
%%%%%%%%%%%%%%%%%%%%%%%%%%%%%%%%%%%%%%%%%%%%%%%%%%%%%%%%%%%%%%%%%%%%%%%%%%%%%%%%%%%%%%%%%%%%%%%%%%%%%%%%%%%%%%%%%%%%%%%%

Let $n\geq 2$ be an integer, and let $D$ be a weakly connected digraph satisfying (D1), (D2) and (D3).
We first suppose that $n=2$.
Since $G_{D}$ is $K_{1,2}$-free, $G_{D}$ is a complete graph, i.e., $\alpha (G_{D})=1$.
Then by Lemma~\ref{lem2.0}, ${\rm pc}(D)={\rm pp}(D)=1$, as desired.
Thus we may assume that $n\geq 3$.

Set $n_{0}=\max\{\lceil \frac{n^{2}-n-2}{2}\rceil ,n\}$.
Take a longest induced path $P$ of $G_{D}$, and write $P=u_{1}u_{2}\cdots u_{m}$.
Let $X_{i}~(i\geq 0)$ and $Y$ be the sets defined as in Subsection~\ref{sec2-1}.
Then by Lemma~\ref{lem2.4}, $V(D)$ is the disjoint union of $V(P)$, $Y$ and $X_{i}~(1\leq i\leq 2n_{0}-1)$.
By Lemmas~\ref{lem2.0} and \ref{lem-alpha},
\begin{align*}
{\rm pc}(D-(V(P)\cup Y)) \leq {\rm pp}(D-(V(P)\cup Y)) &\leq \sum _{1\leq i\leq 2n_{0}-1}{\rm pp}(D[X_{i}])\\
&\leq \sum _{1\leq i\leq 2n_{0}-1}\alpha (G_{D}[X_{i}])\\
&\leq \sum _{1\leq i\leq 2n_{0}-1}\alpha _{i}
\end{align*}
where $\alpha _{i}$ is the constant defined in Subsection~\ref{sec2-1}.
As we mentioned before Lemma~\ref{lem-alpha}, the value $\alpha _{i}$ depends on $n$ and $i$.
Furthermore, $n_{0}$ is a constant depending on $n$ only.
Consequently, the value $\sum _{1\leq i\leq 2n_{0}-1}\alpha _{i}$ depends on $n$.
Thus it suffices to show that 
\begin{enumerate}[{$\bullet $}]
\item
${\rm pc}(D[V(P)\cup Y])$ is bounded by a constant depending on $n$ only, and
\item
if $D$ satisfies (D'1), then ${\rm pp}(D[V(P)\cup Y])$ is bounded by a constant depending on $n$ only.
\end{enumerate}

Let $Q=D[V(P)]$.
Note that $Q$ is a pseudo-path in $D$ and $V(Q)=V(P)$.
Let $I=\{i:1\leq i\leq m,~d^{+}_{Q}(u_{i})=2\mbox{ or }d^{-}_{Q}(u_{i})=2\}$, and write $I=\{i_{1},\ldots ,i_{p}\}$ where $i_{1}<i_{2}<\cdots <i_{p}$.
Then $I\cap \{1,m\}=\emptyset $ and it follows from (D3) that $p=|I|=r(Q)\leq n$.
For each $h$ with $0\leq h\leq p$, let $Q_{h}=Q[\{u_{j}:i_{h}\leq j\leq i_{h+1}\}]$ where $i_{0}=1$ and $i_{p+1}=m$.
Note that $Q_{0}=Q$ if $I=\emptyset $.
By the definition of $I$, $Q_{0},Q_{1},\ldots ,Q_{p}$ are induced directed paths in $D$.

Let
\begin{align*}
J_{0} &= \{j:1\leq j\leq n_{0}\mbox{ or }m-n_{0}+1\leq j\leq m\},\\
J_{1} &= \{j:1\leq h\leq p,~i_{h}-n+1\leq j\leq i_{h}+n-1\},\mbox{ and}\\
Y^{\sharp } &= Y\cap \left(\bigcup _{j\in J_{1}}N_{G_{D}}(u_{j})\right).
\end{align*}
Now we list two facts which will be used in the proof.
\begin{enumerate}
\item[{\bf (Q1)}]
We have $X_{0}=\{u_{i}:i\in J_{0}\}$.
In particular, $Y\subseteq \bigcup _{\substack{1\leq j\leq m\\ j\notin J_{0}}}N_{G_{D}}(u_{j})$.
\item[{\bf (Q2)}]
For $h$ with $0\leq h\leq p$, $V(Q_{h})\setminus \{u_{j}:j\in J_{0}\cup J_{1}\}\subseteq \{u_{j}:i_{h}+n\leq j\leq i_{h+1}-n\}$.
In particular, if $i_{h+1}-i_{h}\leq 2n-1$, then $V(Q_{h})\setminus \{u_{j}:j\in J_{0}\cup J_{1}\}=\emptyset $.
\end{enumerate}
Since $G_{D}$ is $K_{1,n}$-free, we have $\alpha (G_{D}[N_{G_{D}}(u_{j})])\leq n-1$ for every $j\in J_{1}$.
This together with Lemma~\ref{lem2.0} implies the following:
\begin{align}
{\rm pc}(D[Y^{\sharp }]) \leq {\rm pp}(D[Y^{\sharp }]) &\leq \alpha (G_{D}[Y^{\sharp }])\nonumber \\
&\leq \alpha \left(G_{D}\left[\bigcup _{j\in J_{1}}N_{G_{D}}(u_{j})\right]\right)\nonumber \\
&\leq \sum _{j\in J_{1}}\alpha (G_{D}[N_{G_{D}}(u_{j})])\nonumber \\
&\leq |J_{1}|(n-1)\nonumber \\
&\leq p(2n-1)(n-1)\nonumber \\
&\leq n(2n-1)(n-1).\label{eq-pf-mainthm-partition}
\end{align}

For each $h$ with $0\leq h\leq p$, let $Z_{h}=N_{G_{D}}(V(Q_{h})\setminus \{u_{j}:j\in J_{0}\cup J_{1}\})\cap (Y\setminus Y^{\sharp })$.
Now we prove that
\begin{align}
\mbox{$Y\setminus Y^{\sharp }$ is the disjoint union of $Z_{0},Z_{1},\ldots ,Z_{p}$.}\label{cond-proof-mainthm-YYsharp}
\end{align}
By (Q1) and the definition of $Y^{\sharp }$, we have $Y\setminus Y^{\sharp }\subseteq \bigcup _{\substack{1\leq j\leq m\\ j\notin J_{0}\cup J_{1}}}N_{G_{D}}(u_{j})$, and hence $Y\setminus Y^{\sharp }=\bigcup _{0\leq h\leq p}Z_{h}$.
Thus it suffices to show that $Z_{0},Z_{1},\ldots ,Z_{p}$ are pairwise disjoint.
By way of contradiction, we suppose that $Z_{h}\cap Z_{h'}\neq \emptyset $ for some $h$ and $h'$ with $0\leq h<h'\leq p$.
Take $z\in Z_{h}\cap Z_{h'}$.
Then there exist two indices $a$ and $b$ such that $zu_{a},zu_{b}\in E(G_{D})$, $i_{h}\leq a\leq i_{h+1}\leq i_{h'}\leq b\leq i_{h'+1}$ and $a,b\notin J_{0}\cup J_{1}$.
In particular, $a<i_{h+1}<b$.
Then $1\leq h+1\leq p$, and so $\{j:i_{h+1}-n+1\leq j\leq i_{h+1}+n-1\}\subseteq J_{1}$.
This implies that $a\leq i_{h+1}-n<i_{h+1}<i_{h+1}+n\leq b$, and hence $b-a\geq (i_{h+1}+n)-(i_{h+1}-n)=2n\geq 6$.
On the other hand, since $z$ is adjacent to no vertices in $X_{0}$, applying Lemma~\ref{lem-2-3-0}(i) with $y=z$ and $Q=P$, we have $b-a\leq 3$, which is a contradiction.
Consequently, (\ref{cond-proof-mainthm-YYsharp}) holds.

Recall that $Q_{h}$ is a directed path in $D$, and $u_{i_{h}}$ is either the internal vertex or terminal vertex of $Q_{h}$.
By (Q2), $Z_{h}$ is a subset of $N_{G_{D}}(V(Q_{h}))$ such that $N_{G_{D}}(z)\cap V(Q_{h})\subseteq \{u_{j}:i_{h}+n\leq j\leq i_{h+1}-n\}$ for every $z\in Z_{h}$.
Thus, applying Lemma~\ref{lem-2-3-1} with $Q=Q_{h}$ and $Y'=Z_{h}$, we obtain the following:
\begin{enumerate}
\item[{\bf (P1)}]
For every $h$ with $0\leq h\leq p$, ${\rm pc}(D[V(Q_{h})\cup Z_{h}])\leq \frac{(n-2)n(n+5)}{2}+6(n-2)$.
\item[{\bf (P2)}]
If $D$ satisfies (D'1), then
\begin{enumerate}[{$\circ$}]
\item
${\rm pp}(D[V(Q_{0})\cup Z_{0}])\leq \frac{n(n+5)}{2}+1$ and
\item
for every $h$ with $1\leq h\leq p$, ${\rm pp}(D[(V(Q_{h})\setminus \{u_{i_{h}}\})\cup Z_{h}])\leq \frac{n(n+5)}{2}+1$.
\end{enumerate}
\end{enumerate}
By (\ref{cond-proof-mainthm-YYsharp}), we can easily verify that $V(P)\cup (Y\setminus Y^{\sharp })$ is the disjoint union of
$$
V(Q_{0})\cup Z_{0},(V(Q_{1})\setminus \{u_{i_{1}}\})\cup Z_{1},\ldots ,(V(Q_{p})\setminus \{u_{i_{p}}\})\cup Z_{p}.
$$
This together with (\ref{eq-pf-mainthm-partition}), (P1) and (P2) implies that
\begin{align*}
{\rm pc}(D[V(P)\cup Y]) &\leq {\rm pc}(D[Y^{\sharp }])+\sum _{0\leq h\leq p}{\rm pc}(D[V(Q_{h})\cup Z_{h}])\\
&\leq n(2n-1)(n-1)+(p+1)\left(\frac{(n-2)n(n+5)}{2}+6(n-2)\right)\\
&\leq n(2n-1)(n-1)+(n+1)\left(\frac{(n-2)n(n+5)}{2}+6(n-2)\right)
\end{align*}
and, if $D$ satisfies (D'1), then
\begin{align*}
{\rm pp}(D[V(P)\cup Y]) &\leq {\rm pp}(D[Y^{\sharp }])+{\rm pp}(D[V(Q_{0})\cup Z_{0}])+\sum _{1\leq h\leq p}{\rm pp}(D[(V(Q_{h})\setminus \{u_{i_{h}}\})\cup Z_{h}])\\
&\leq n(2n-1)(n-1)+(p+1)\left(\frac{n(n+5)}{2}+1\right)\\
&\leq n(2n-1)(n-1)+(n+1)\left(\frac{n(n+5)}{2}+1\right).
\end{align*}
This completes the proof of Theorem~\ref{mainthm}.

%%%%%%%%%%%%%%%%%%%%%%%%%%%%%%%%%%%%%%%%%%%%%%%%%%%%%%%%%%%%%%%%%%%%%%%%%%%%%%%%%%%%%%%%%%%%%%%%%%%%%%%%%%%%%%%%%%%%%%%%
%%%%%%%%%%%%%%%%%%%%%%%%%%%%%%%%%%%%%%%%%%%%%%%%%%%%%%%%%%%%%%%%%%%%%%%%%%%%%%%%%%%%%%%%%%%%%%%%%%%%%%%%%%%%%%%%%%%%%%%%
%%%%%%%%%%%%%%%%%%%%%%%%%%%%%%%%%%%%%%%%%%%%%%%%%%%%%%%%%%%%%%%%%%%%%%%%%%%%%%%%%%%%%%%%%%%%%%%%%%%%%%%%%%%%%%%%%%%%%%%%
\section{Ramsey-type theorem for the cycle cover/partition number of digraphs}\label{sec4}
%%%%%%%%%%%%%%%%%%%%%%%%%%%%%%%%%%%%%%%%%%%%%%%%%%%%%%%%%%%%%%%%%%%%%%%%%%%%%%%%%%%%%%%%%%%%%%%%%%%%%%%%%%%%%%%%%%%%%%%%
%%%%%%%%%%%%%%%%%%%%%%%%%%%%%%%%%%%%%%%%%%%%%%%%%%%%%%%%%%%%%%%%%%%%%%%%%%%%%%%%%%%%%%%%%%%%%%%%%%%%%%%%%%%%%%%%%%%%%%%%
%%%%%%%%%%%%%%%%%%%%%%%%%%%%%%%%%%%%%%%%%%%%%%%%%%%%%%%%%%%%%%%%%%%%%%%%%%%%%%%%%%%%%%%%%%%%%%%%%%%%%%%%%%%%%%%%%%%%%%%%

A digraph $C$ is called a {\it directed cycle} if $G_{C}$ is a cycle and $d^{+}_{C}(x)=d^{-}_{C}(x)=1$ for every $x\in V(C)$.
Let $\Gamma $ be a (di)graph.
When $\Gamma $ is a graph, a family $\CC$ of subgraphs of $\Gamma $ is called a {\it cycle cover} of $\Gamma $ if $\bigcup _{C\in \CC}V(C)=V(\Gamma )$ and each element of $\CC$ is either a cycle or a connected graph of order at most two.
When $\Gamma $ is a digraph, a family $\CC$ of subdigraphs of $\Gamma $ is called a {\it cycle cover} of $\Gamma $ if $\bigcup _{C\in \CC}V(C)=V(\Gamma )$ and each element of $\CC$ is either a directed cycle or a weakly connected digraph of order at most two.
A cycle cover $\CC$ of $\Gamma $ is called a {\it cycle partition} of $\Gamma $ if the elements of $\CC$ are pairwise vertex-disjoint.
The value $\min \{|\CC |:\CC\mbox{ is a cycle cover of }\Gamma \}$ (resp. $\min \{|\CC |:\CC\mbox{ is a cycle partition of }\Gamma \}$), denoted by ${\rm cc}(\Gamma )$ (resp. ${\rm cp}(\Gamma )$), is called the {\it cycle cover number} (resp. the {\it cycle partition number}) of $\Gamma $.
In \cite{CF}, the authors proved the following Ramsey-type result for the cycle partition number of graphs.
(They also showed that the $\{K^{*}_{n},K_{1,n},P_{n}\}$-freeness in the theorem is necessary to bound the cycle cover number by a constant.
Thus Theorem~\ref{ThmB} also offers tight forbidden subgraph conditions on a Ramsey-type theorem for the cycle cover number of graphs.)

\begin{Thm}[Chiba and Furuya~\cite{CF}]%%%%%%%%%%%%%%%%%%%%%%%%%%%%%%%%%%%%%%%%%%%%%%%%%%%%%%%%%%%%%%%%%%%%%%%%%%%%%%%%%
\label{ThmB}
For an integer $n\geq 2$, there exists a constant $c_{0}=c_{0}(n)$ depending on $n$ only such that ${\rm cp}(G)\leq c_{0}$ for every connected $\{K^{*}_{n},K_{1,n},P_{n}\}$-free graph $G$.
\end{Thm}
%%%%%%%%%%%%%%%%%%%%%%%%%%%%%%%%%%%%%%%%%%%%%%%%%%%%%%%%%%%%%%%%%%%%%%%%%%%%%%%%%%%%%%%%%%%%%%%%%%%%%%%%%%%%%%%%%%%%%%%%

In this section, we give an analogy of Theorem~\ref{ThmB} for digraphs.
Let $T_{n}$ denote a {\it transitive tournament} of order $n$, i.e., $T_{n}$ is a digraph on $\{1,2,\ldots ,n\}$ such that $E(T_{n})=\{(i,j):1\leq i<j\leq n\}$.
Our main result in this section is the following.

\begin{thm}%%%%%%%%%%%%%%%%%%%%%%%%%%%%%%%%%%%%%%%%%%%%%%%%%%%%%%%%%%%%%%%%%%%%%%%%%%%%%%%%%%%%%%%%%%%%%%%%%%%%%%%%%%%%%
\label{mainthm-2}
For an integer $n\geq 2$, there exists a constant $c_{1}=c_{1}(n)$ depending on $n$ only such that ${\rm cp}(D)\leq c_{1}$ for every weakly connected $T_{n}$-free digraph $D$ such that $G_{D}$ is $\{K_{1,n},P_{n}\}$-free.
\end{thm}
%%%%%%%%%%%%%%%%%%%%%%%%%%%%%%%%%%%%%%%%%%%%%%%%%%%%%%%%%%%%%%%%%%%%%%%%%%%%%%%%%%%%%%%%%%%%%%%%%%%%%%%%%%%%%%%%%%%%%%%%

Since $T_{n}$ has no directed cycles, we have ${\rm cp}(T_{n})=\lceil \frac{n}{2} \rceil $.
In particular, when we consider a forbidden subdigraph condition assuring us that the cycle cover/partition number of digraphs is bounded by a constant, some transitive tournaments appear in the condition.
Thus the forbidden structure condition in Theorem~\ref{mainthm-2} is best possible in a sense.

\medbreak\noindent\textit{Proof of Theorem~\ref{mainthm-2}.}\quad
Let $D$ be a weakly connected $T_{n}$-free digraph such that $G_{D}$ is $\{K_{1,n},P_{n}\}$-free.
If $G_{D}$ contains a clique $C$ with $|C|=2^{n-1}$, then it is well-known that $D$ contains a copy of $T_{n}$ as an induced subdigraph, which is a contradiction.
Thus $G_{D}$ is $K_{2^{n-1}}$-free.
For every connected graph $G$, it is known that $|V(G)|\leq \Delta (G)^{\mbox{diam}(G)}$, where $\Delta (G)$ and ${\rm diam}(G)$ are the maximum degree and the diameter of $G$, respectively.
Since $G_{D}$ is $\{K_{2^{n-1}},K_{1,n}\}$-free, we can easily verify that the maximum degree of $G_{D}$ is at most $R(2^{n-1}-1,n)-1$.
Since $G_{D}$ is $P_{n}$-free, the diameter of $G_{D}$ is at most $n-2$.
Consequently, ${\rm cp}(D)\leq |V(D)|=|V(G_{D})|\leq \Delta (G_{D})^{\mbox{diam}(G_{D})}\leq (R(2^{n-1}-1,n)-1)^{n-2}$, as desired.
\qed

One might notice that the $K^{*}_{n}$-freeness appearing in Theorem~\ref{ThmB} is absorbed by the $T_{n}$-freeness in Theorem~\ref{mainthm-2}.
Thus the forbidden structure in Theorem~\ref{mainthm-2} does not form the forbidden subgraph condition in Theorem~\ref{ThmB} plus additional assumptions which are peculiar to digraphs.

%%%%%%%%%%%%%%%%%%%%%%%%%%%%%%%%%%%%%%%%%%%%%%%%%%%%%%%%%%%%%%%%%%%%%%%%%%%%%%%%%%%%%%%%%%%%%%%%%%%%%%%%%%%%%%%%%%%%%%%%
%%%%%%%%%%%%%%%%%%%%%%%%%%%%%%%%%%%%%%%%%%%%%%%%%%%%%%%%%%%%%%%%%%%%%%%%%%%%%%%%%%%%%%%%%%%%%%%%%%%%%%%%%%%%%%%%%%%%%%%%
%%%%%%%%%%%%%%%%%%%%%%%%%%%%%%%%%%%%%%%%%%%%%%%%%%%%%%%%%%%%%%%%%%%%%%%%%%%%%%%%%%%%%%%%%%%%%%%%%%%%%%%%%%%%%%%%%%%%%%%%
\section*{Acknowledgment}
%%%%%%%%%%%%%%%%%%%%%%%%%%%%%%%%%%%%%%%%%%%%%%%%%%%%%%%%%%%%%%%%%%%%%%%%%%%%%%%%%%%%%%%%%%%%%%%%%%%%%%%%%%%%%%%%%%%%%%%%
%%%%%%%%%%%%%%%%%%%%%%%%%%%%%%%%%%%%%%%%%%%%%%%%%%%%%%%%%%%%%%%%%%%%%%%%%%%%%%%%%%%%%%%%%%%%%%%%%%%%%%%%%%%%%%%%%%%%%%%%
%%%%%%%%%%%%%%%%%%%%%%%%%%%%%%%%%%%%%%%%%%%%%%%%%%%%%%%%%%%%%%%%%%%%%%%%%%%%%%%%%%%%%%%%%%%%%%%%%%%%%%%%%%%%%%%%%%%%%%%%

This work was partially supported by JSPS KAKENHI Grant number JP20K03720 (to S.C) and JSPS KAKENHI Grant number JP18K13449 (to M.F).

\end{document}